\documentclass{article}
\usepackage[english]{babel}
\usepackage[utf8]{inputenc}

\usepackage{bm}
\usepackage{amsmath}
\usepackage{graphicx}
\usepackage{color}
\usepackage{multirow}
\usepackage{amsfonts}
\newcommand{\rv}[1]{\textcolor{black}{#1}}
\newcommand{\grad}{\mathop{\rm grad}\nolimits}
\renewcommand{\div}{\mathop{\rm div}\nolimits}

\def\bu{\mathbf{u}}
\def\bv{\mathbf{v}}
\def\bsig{\boldsymbol{\sigma}}
\def\bvar{\boldsymbol{\varepsilon}}
\graphicspath{ {./figs/} } 

\date{}

\author{
Aleksei Tyrylgin
\thanks{Laboratory of Computational Technologies for Modeling Multiphysical and Multiscale Permafrost Processes, North-Eastern Federal University, Yakutsk, Republic of Sakha (Yakutia), 677980, Russia \& North-Caucasus Center for Mathematical Research, North-Caucasus Federal University,
Stavropol, 355017, Russia.}
\and
Maria Vasilyeva \thanks{Department of Mathematics and Statistics, Texas A\&M University, Corpus Christi, 78412, Texas, USA. Email: {\tt maria.vasilyeva@tamucc.edu}.}
\and
Anatoly Alikhanov \thanks{North-Caucasus Center for Mathematical Research, North-Caucasus Federal University,
Stavropol, 355017, Russia.  Email: {\tt aaalikhanov@gmail.com}.}
\and
Dongwoo Sheen \thanks{Department of Mathematics, Seoul National
  University, Seoul 08826, South Korea.  Email: {\tt dongwoosheen@gmail.com}.}
}

\begin{document}
\title{A computational macroscale model for the time fractional poroelasticity problem in fractured and heterogeneous media}

\maketitle

\begin{abstract}
In this work, we introduce a time memory formalism in poroelasticity model that couples
the pressure and displacement. We assume this multiphysics process occurs in multicontinuum media. 
The mathematical model contains a coupled system of equations for pressures in each continuum and elasticity equations for displacements of the medium. We assume that the temporal dynamics is governed by fractional derivatives following some works in the literature. We derive an implicit  finite difference approximation 
for time discretization based on the  Caputo’s time fractional  derivative. A Discrete Fracture Model (DFM) is used to model fluid flow through fractures and treat the complex network of fractures. 
We assume different fractional powers in fractures and matrix due to slow and fast dynamics. We develop a coarse grid approximation based on the Generalized Multiscale Finite Element Method (GMsFEM), where we solve local spectral problems for construction of the multiscale basis functions. We present numerical results  for the two-dimensional model problems in fractured heterogeneous porous media. We investigate error analysis between reference (fine-scale) solution and multiscale solution with different numbers of multiscale basis functions.
 The results show that the proposed method can provide good accuracy on a coarse grid. 
\end{abstract}

\section{Introduction}

Recently, various applications of \rv{differential equations with fractional order derivatives have been of increasing interest.} Moreover, in contrast to the classical derivative of integer order, there are many non-identical definitions of derivatives of fractional order \cite{du2020temporal, gu2017fast, gao2017temporal}. For example, fractional derivative is used in models of viscoelastic bodies, continuous media, transformation of temperature, humidity in atmospheric layers, diffusion equations, and in other areas \cite{liu2015he, mahiuddin2020application, zaky2021numerical, imran2017heat, vasil2021iterative}. \rv{In addition, in the middle of the twentieth century F. Mainardi and M. Caputo showed that the use of differential equations with fractional derivatives more better the models of thermoviscoelasticity are described, which make it possible to more accurately reproduce the experimentally observed data \cite{caputo1971linear, caputo1974vibrations, mondal2020interactions, shariyat20203d}.} The relevance of such a study is explained by the fact that the use of a rich arsenal of fractional differentiation methods will make new methods for statistical analysis of nonstationary time series. 

In this paper, we study \rv{a} fractional poroelasticity model. The basic
mathematical structure of the poroelasticity models are usually
coupled equations for pressure and displacement. In their modern form,
such models were proposed in the works of M. Biot
\cite{biot1941general, biot1955theory, castelletto2018multiscale,
  gaspar2008stabilized, kolesov2014splitting}. Biot's \rv{model
  describes} the coupled processes of deformation of the elastic
medium and the \rv{flow of fluid}. The model is macroscopic in the
sense that the space containing the poroelastic medium is filled with
a two-phase medium,
with one phase corresponding directly to the porous medium,
and the second to the fluid contained in the pores \cite{meirmanov2014mathematical, iliev2016numerical}. In our case, we consider mathematical models for flow in the multicontinuum media, which describes complex flow processes in multiscale fractured heterogeneous porous media \cite{d2008coupling, d2012mixed, dpgmsfem2017, akkutlu2018multiscale, akkutlu2015multiscale}. The flow in fractures has a significant impact on filtration processes and requires careful consideration \cite{lee2001hierarchical, li2006efficient}. \rv{Since }the fractures are characterized by high permeability and their thickness is significantly smaller than the size of the simulated field, this leads to the need to build special mathematical models of multicontinuum, where independent variables are distinguished to describe the flow in a porous medium and in the network of fractures taking into account the special flow function \cite{akkutlu2018multiscale, karimi2004efficient, li2008efficient}. 

The extension of the poroelasticity to include fractional time derivatives appeared in earlier works
\cite{alaimo2019fractional,carcione2010computational,lorenzi2014direct,enelund2006time}, where physical motivations are presented. The fractional time derivatives represent the memory effects
that occur in porous media flows. One way to account for them is to introduce fractional time derivatives.
The main distinction of our model from previous models consists of
several points. First, we use different fractional time derivatives
for each continua as each continua can have different propagation
dynamics, \rv{and thus,} different memory terms. Secondly, we assume that the media have
multiple spatial scales, which commonly occur in porous media
applications. Our main \rv{goal is to derive} a computational
macroscopic model for
fractional Biot's system in multicontinuum heterogeneous media.

Our computational macroscale model is based on \rv{the} Generalized Multiscale Finite Element method (GMsFEM).
\rv{The} GMsFEM has been studied for a various applications related to
poroelasticity problems \cite{brown2016generalized1,
  brown2016generalized2, akkutlu2018multiscale,
  tyrylgin2020generalized}. \rv{The} multiscale finite volume method
\rv{has been applied for the simulation of} the flow problems in fractured
porous media \cite{hajibeygi2011hierarchical, ctene2016algebraic}. For the
effective numerical solution of such problems different homogenization
techniques \rv{have been} developed \cite{li2008efficient, li2006efficient,
  arbogast1990derivation, talonov2016numerical,
  bakhvalov1984homogenization, gavrilieva2018numerical}. Mathematical
models of the flow problems in fractured porous media using \rv{the}
GMsFEM have been researched \cite{yao2010discrete,
  vasilyeva2019nlmcdp, tene2015algebraic,
  dpgmsfem2017}. \rv{The} GMsFEM and NLMC approach for
solution of the flow problems in multicontinuum media \rv{have been
  generalized} in \cite{vasilyeva2019nlmcdp, vasilyeva2019nonlocal, poroelnlmc2018}.

In this paper, we consider \rv{the} GMsFEM for the poroelasticity 
problems in multicontinuum media with the fractional  order time derivative.
For temporal discretization, we use a finite different approximation,
which has a memory term. Since our model equation has multiple fractional
powers, there multiple unknowns with memory effects. 
Because the media properties have multiple scales, we 
use multiscale basis functions to reduce the dimension or the problem.
The multiscale basis functions	are constructed for flow and mechanics. 
Construction of the basis functions for flow problem in multicontinuum 
media is based on the solution of the
coupled system of equations in each local domains. In each coarse grid block, 
we compute the snapshots by solving local problems for pressures in multicontinuum media and displacement. Taking the corresponding to the dominant eigenvalues, and multiplying by a multiscale partition of unity, \rv{we can construct} our multiscale basis functions.

Numerical results are presented. In our numerical examples, we consider
two different type of media. In the first case, we have one continuum.
In the second example, we have two continuum, which increases the number 
of equations. The memory term is handled by saving solution snapshots.
Because the solution is on the coarse grid, this saves some computational
time. We consider numerical simulations using different number of basis
functions. In all examples, when we increase the number of basis functions,
the error decreases. In particular, using fewer basis functions, we obtain 
accurate solution approximation.

The work is organized as follows. In Section 2, we present the
mathematical model of the poroelasticity problem in multicontinuum
medium. Then in Section 3, a fine grid approximation is constructed
using the finite element method and the fractional-order
derivative. In Section 4, we present a coarse grid approximation using
\rv{the} Generalized Multiscale Finite Element method, where we
describe \rv{the} construction of the multiscale basis functions and coarse grid system construction. Numerical results for two-dimensional model poroelasticity problems are presented in Section 5. Finally, we present conclusions.

\section{Problem formulation}
The time fractional flow in porous media  $\Omega \subset R^d$ can be described by the following equation
\begin{equation}
\label{m}
c \frac{\partial^{\alpha} p}{\partial t^{\alpha}}
- \nabla \cdot (k \nabla p) = 0,  \quad \rv{\Omega\times(0,T)}, 
\end{equation}
where 
\rv{$k = \frac{\kappa}{\mu}$ with the permeability $k$ and 
the fluid viscosity $\mu$}, and 
\begin{equation}
\partial^\alpha_t p\rv{(t)} = \frac{1}{\Gamma(1-\alpha)}\int_0^t (t-s)^{-\alpha} \frac{\partial p}{\partial s}\rv{(s)} ds, \quad 0 < \alpha \leq 1,
\end{equation}
\rv{denotes} the Caputo derivative of the order $\alpha$. 

To consider flow in fractured porous media, we \rv{denote $\gamma \subset R^{d-1}$} as a computational domain for low dimensional fracture networks model.  
Therefore, we have the following system of equations for flow in fractured porous media:
\begin{equation}
\label{mf}
\begin{split}
c_m \frac{\partial^{\alpha_m} p_m}{\partial t^{\alpha_m}}
- \nabla \cdot (k_m \nabla p_m) + \eta_{mf} (p_m - p_f)= 0,  
\quad \rv{\Omega\times(0,T)}, \\
c_f \frac{\partial^{\alpha_f} p_f}{\partial t^{\alpha_f}}
- \nabla \cdot (k_f \nabla p_f)  + \eta_{mf} (p_f - p_m) = 0,
 \quad \rv{\gamma\times(0,T)},
\end{split}
\end{equation}
where 
$p_m$ and $p_f$ are the pressure in porous matrix and fractures, 
$\kappa_m$ and $\kappa_f$ the porous matrix and fractures permeability ($k_m = \frac{\kappa_m}{\mu}$ and $k_f = \frac{\kappa_f}{\mu}$), \rv{and}
$\eta_{mf}$ the mass transfer term between the porous matrix and fractures. 

We can write \rv{a} similar system of equations for flow in triple continuum approach,  where the first continuum  describe a flow in the matrix of the porous media,  the second continuum belongs to the network of small highly connected fracture network (natural fractures) and the third continuum related to the flow in low dimensional fracture networks (large-scale fractures). 
We have following system of equations for $(p_1, p_2, p_f)$: 
\begin{equation}
\label{mmf}
\begin{split}
c_1  \frac{\partial^{\alpha_1} p_1}{\partial t^{\alpha_1}}
-  \nabla \cdot (k_1 \nabla p_1)
+ \eta_{12} (p_1 - p_2) 
+ \eta_{1f} (p_1 - p_f) = 0,  \quad \rv{\Omega\times(0,T)}, \\
c_2  \frac{\partial^{\alpha_2} p_2}{\partial t^{\alpha_2}}
-  \nabla \cdot (k_2 \nabla p_2) 
+ \eta_{12} (p_2 - p_1)
+ \eta_{2f} (p_2 - p_f) = 0, \quad \rv{\Omega\times(0,T)},
\\
c_f \frac{\partial^{\alpha_f} p_f}{\partial t^{\alpha_f}}
-  \nabla \cdot (k_f \nabla p_f)  
+ \eta_{1f} (p_f - p_1)
+ \eta_{2f} (p_f - p_2) = 0,  \quad \rv{\gamma\times(0,T)},
\end{split}
\end{equation}
where, \rv{for the continuum index $i = 1,2,f,$}
$p_i$ denotes the pressure, 
$\kappa_i$ the permeability ($k_i = \frac{\kappa_i}{\mu}$, $\mu$ the fluid viscosity), and
$\eta_{ij}$ the mass transfer term that are proportional to the continuum permeabilities.

We can generalize it as flow model for multicontinuum media
\begin{equation}
\label{mc}
c_i \frac{\partial^{\alpha_i} p_i}{\partial t^{\alpha_i}} 
 - \nabla \cdot (k_i \nabla p_i)  + \sum_{j \neq i}  \eta_{ij} (p_i - p_j)  = 0, 
\quad \rv{\Omega\times(0,T)},
\end{equation}
where $i = 1,...,M$ and $M$ is the number of continua. 

For the mechanics of the poroelastic multicontinuum media, we use an effective equation for displacement and have following poroelasticity problem for multicontinuum media
\begin{equation}
\label{mcu}
\begin{split}
c_i \frac{\partial^{\alpha_i} p_i}{\partial t^{\alpha_i}} 
+  \gamma_i \frac{\partial^{\beta_i} \div  \bu}{\partial t^{\beta_i}} 
- \nabla \cdot (k_i \nabla p_i)  + \sum_{j \neq i}  \eta_{ij} (p_i - p_j)  &= 0, 
\quad \rv{\Omega\times(0,T)},  \quad i = 1,...,M,\\
 - \nabla  \bsig (\bu) + \sum_{j} \gamma_j \nabla p_j  &= 0 , 
\quad \rv{\Omega\times(0,T)},  
\end{split}
\end{equation}
where  $\bsig$ \rv{denotes} the stress tensor, $\bu$ the displacement,
$\gamma_i$ the Biot coefficient, $M_i$  the Biot modulus ($c_i =
\frac1{M_i}$)
for the $i$-th component. 
In the case of a linear elastic stress-strain constitutive relation, we have
\[
\bsig(\bu) = 2 \mu \bvar(\bu) + \lambda \div \bu \ \mathcal{I}, 
\quad  
\bvar(\bu) = \frac{1}{2}(\nabla \bu + \nabla \bu^T),
\]
where $\bvar$ is the strain tensor, $\lambda$ and $\mu$ are the Lame’s coefficients. 
Here we have a volume force sources that proportional to the sum of the pressure gradients for each continuum.  

In the presented poroelasticity model \eqref{mcu},  the fractional time parameters $\alpha_i$ and $\beta_i$  are used to simulate the effects of history on porous media flow, where $\alpha_i$ is used for the effects of flow (pressure) history and $\beta_i$ for the effects of mechanics (displacements) history on flow processes in  multicontinuum media.  

Next, we will \rv{concentrate on the} triple continuum poroelasticity model:
\begin{equation}
\label{eq:main}
\begin{split}
c_1 \frac{\partial^{\alpha_1} p_1}{\partial t^{\alpha_1}} 
+ \gamma_1 \frac{\partial^{\beta_1} \div  \bu}{\partial t^{\beta_1}} 
- \nabla \cdot (k_1 \nabla p_1)  
+ \eta_{12} (p_1 - p_2) + \eta_{1f} (p_1 - p_f) &= 0, 
\quad \rv{\Omega\times(0,T)}, \\
c_2 \frac{\partial^{\alpha_2} p_2}{\partial t^{\alpha_2}} 
+ \gamma_2 \frac{\partial^{\beta_2} \div  \bu}{\partial t^{\beta_2}}
- \nabla \cdot (k_2 \nabla p_2)  
+ \eta_{12} (p_2 - p_1) + \eta_{2f} (p_2 - p_f) &= 0, 
\quad \rv{\Omega\times(0,T)}, \\ 
c_f \frac{\partial^{\alpha_f} p_f}{\partial t^{\alpha_f}} 
+ \gamma_f \frac{\partial^{\beta_f} \div  \bu}{\partial t^{\beta_f}}
- \nabla \cdot (k_f \nabla p_f)  
+ \eta_{1f} (p_f - p_1) + \eta_{2f} (p_f - p_2) &= 0, 
 \quad \rv{\gamma\times(0,T)}, \\
 - \div \bsig (\bu) 
 + \gamma_1 \nabla p_1  + \gamma_2  \nabla p_2 + \gamma_f \nabla p_f  &= 0 , 
\quad \rv{\Omega\times(0,T)}, 
\end{split}
\end{equation}
where the first continuum \rv{describes} a flow in the matrix of the porous media, the second continuum belongs to the network of small highly connected fracture network (natural fractures) and the third continuum \rv{relates} to the flow in low dimensional fracture networks.

\section{Fine grid approximation using FEM}
For \rv{the temporal} approximation, we use an uniform mesh with $N_T$ time steps and time step size $\tau = \frac{T}{N_T}$, where $T$ is the final time for simulation.  The values of a pressures and displacement at the time $t^n = n \tau$ ($n=0,1,2,..,N_T$) are denoted by $(p_i^n, \bu^n) = (p_i(t^n), \bu(t^n))$, where $p_i$ is the pressure of the $i$-th continuum. 

The fractional-order derivative of the fucntion $v^n$ is defined using
the following formula: \rv{add references here}
\[
\frac{\partial^{\alpha} v^n}{\partial t^{\alpha}}
\approx \zeta^{(\alpha)}_{\tau} \left( 
v^n - v^{n-1} + \sum_{j=2}^n \zeta^{(\alpha)}_{j-1} (v^{n-j+1} - v^{n-j})
\right),
\]
where 
\[
\zeta^{(\alpha)}_{\tau} = \frac{1}{\tau^{\alpha} \Gamma(2 - \alpha)}, \quad
\zeta^{(\alpha)}_{j-1} = j^{1-\alpha} - (j-1)^{1-\alpha}.
\]

For \rv{the spatial} approximation, we use \rv{the} finite element method.  
Let  $V = [H^1(\Omega)]^d$,  
$W_1 = W_2 =  H^1(\Omega)$ and 
$W_f = H^1(\gamma)$. 
The variational formulation of  
the poroelasticity problem in multicontinuum media \eqref{eq:main} can
be written as follows: \rv{
given $(p_1^0, p_2^0, p_f^0, \bu^0) \in W_1 \times W_2 \times W_f\times V$ iteratively}
find $(p_1^n, p_2^n, p_f^n, \bu^n) \in W_1 \times W_2 \times W_f \times V$ such that
\begin{equation}
\label{eq:fine}
\begin{split}
\zeta^{(\alpha_i)}_{\tau} 
m_i ( p^n_i - p^{n-1}_i, w_i) 
& + \zeta^{(\alpha_i)}_{\tau} 
\sum_{j=2}^n \zeta^{(\alpha_i)}_{j-1} \ m_i(  p^{n-j+1}_i -  p^{n-j}_i, w_i ) \\
& +
\zeta^{(\beta_i)}_{\tau}
d_i(\bu^n -  \bu^{n-1}, w_i ) 
+ \zeta^{(\beta_i)}_{\tau}
\sum_{j=2}^n \zeta^{(\beta_i)}_{j-1} d_i( \bu^{n-j+1} -  \bu^{n-j}, w_i) \\
&+ b_i(p_i^n, w_i) 
+ \sum_{j \neq i} q_{ij}(p_i^n - p_j^n, w_i) =  0,  
\quad \forall w_i \in W_i, \rv{\, i=1,2,f},
\\
& a(\bu^n, \bv) + \sum_j g_i(p_i^n, \bv) = 0, 
\quad \forall \bv \in V,  
\end{split}
\end{equation}
where
\[
b_i(p_i, w_i) = \int_{\Omega_i} k_i \nabla p_i \cdot \nabla w_i dx,  \quad
a(\bu,\bv) = \int_{\Omega} {\bsig}(\bu) \cdot {\bvar}(\bv) \, dx, 
\]\[
m_i(p_i, w_i) =  \int_{\Omega_i} c_i p_i w_i \, dx, \quad 
q_{ij}(p_i - p_j, w_i) =  \int_{\Omega_i} \eta_{ij}(p_i - p_j) \, w_i \, dx,
\]\[
d_i(\bu, w_i) =  \int_{\Omega_i} \gamma_i \div \bu \, w_i \, dx,\quad 
g_i(p_i, \bv) =  \int_{\Omega_i} \gamma_i \nabla p_i \bv \, dx,
\]
\rv{for $i,j=1,2,f,$}
with $\Omega_1 = \Omega_2 = \Omega$, $\Omega_f = \gamma.$

Let $\mathcal{T}^h$ denote a finite element partition of the domain $\Omega$ and $\mathcal{E}_h$ is the set of all the interfaces between the elements $\mathcal{T}_h$. 
For the fracture continuum, we use a discrete fracture model and use an unstructured fine grid $\mathcal{T}^h$ that explicitly resolve fracture geometry.
We assume that $\mathcal{E}_{\gamma} = \cup_j \gamma_j$ is the subset of faces for $\mathcal{T}^h$ that represent fractures, where $j = 1,...,N_{frac}$, $N_{frac}$ is the number of discrete fractures and $\mathcal{E}_{\gamma} \subset \mathcal{E}_h$ be the subset of all faces that represent fractures. Moreover, $\mathcal{E}_{\gamma}$ describe the lower dimensional fracture grid.  

\rv{For $i=1,2,f$, let}
\[
p_i = \sum_l p_{i,l}^h \phi^i_l, \quad  
\bu = \sum_l u_l^h \Phi_l,
\] 
where 
\rv{$\{\Phi_l\}$ is} the basis for displacements,
\rv{$\{\phi^i_l\},i=1,2,$} the $d$-dimensional bases for pressure,
and \rv{$\{\phi^f_l\}$  the $(d-1)$-dimensional basis for pressure.}
Then we have following discrete system in matrix form on the fine grid for the triple-continuum media 
\begin{equation}
\begin{split}
\zeta^{(\alpha_i)}_{\tau} M_i p^n_i + 
& \zeta^{(\beta_i)}_{\tau} D_i \bu^n + 
A_i p^n_i + 
\sum_{j\rv{\neq} i}Q_{ij}(p^n_i - p^n_j) \\
&= 
\zeta^{(\alpha_i)}_{\tau} M_i p^{n-1}_i -
\zeta^{(\alpha_i)}_{\tau} 
\sum_{j=2}^{n} \zeta^{(\alpha_i)}_{j-1} M_i (p^{n-j+1}_i - p^{n-j}_i) \\
& + 
\zeta^{(\beta_i)}_{\tau} D_i \bu^{n-1} -
\zeta^{(\beta_i)}_{\tau}
\sum_{j=2}^{n} \zeta^{(\beta_i)}_{j-1} D_i (\bu^{n-j+1} - \bu^{n-j}),
\rv{\quad\text{ for }i=1,2,f,}
\\
&\sum_{j} D_j^T p^n_j + A_u \bu^n = 0,
\end{split}
\end{equation}
where \rv{the indices $i$ and $n$ stand for the continuum and the time
step, respectively, and}
\[ 
A_i = [a_{i,ln}], \quad 
a_{i, ln} = \int_{\Omega_i} k_i \nabla \phi^i_l \cdot \nabla \phi^i_n dx, 
\quad
A_u = [a_{u,ln}], \quad  
a_{u,ln} = \int_{\Omega} {\bsig}(\Phi_l) \cdot {\bvar}(\Phi_n) dx,
\]\[
M_i=[c_{i,ln}], \quad 
m_{i,ln}=\int_{\Omega_i} c_i \phi^i_l \phi^i_n dx,  
\quad 
Q_{ij} = [q_{ij, ln}], \quad 
q_{ij, ln} = \int_{\Omega_i} \eta_{ij} \phi^i_l \phi^j_n dx, 
\]\[ 
D_i = [d_{i,ln}], \quad 
d_{i,ln} = \int_{\Omega} \gamma_i  \div \Phi_l \phi^i_n  dx. 
\]


\section{Coarse grid approximation using GMsFEM}
For the coarse grid approximation, we use \rv{the} Generalized Multiscale Finite Element Method (GMsFEM).
We construct multiscale basis functions for displacements and
pressures separately, but basis functions for multicontinuum pressure
equations are constructed in \rv{a} coupled way. 

\rv{Denote by $\mathcal{T}^H$} the coarse grid partitioning of the domain
\[
\mathcal{T}^H = \bigcup_j K_j,
\]
where $K_j$'s are coarse grid cells. We will use \rv{the standard
  continuous $P_1$} Galerkin approximation on the coarse grid, and define local domain $\omega_l$ for multiscale basis functions as combination of the several coarse grid cells that share same coarse grid nodes ($l = 1,...,N_v^H$, $N_v^H$ \rv{being} the number of coarse grid vertices).

\subsection{Multiscale basis functions for pressures in multicontinuum media}
To construct a snapshot space, we solve \rv{the} following local problem in domain $\omega_l$:
find \rv{$\psi^{l, j} = (\psi^{l, j}_1,\psi^{l, j}_2,\psi^{l, j}_f)$
  $\in W^h_1 \times W^h_2\times W^h_f$} such that
\begin{equation}
b_i(\psi^{l, j}_i, w_i) + \sum_j q_{ij}(\psi^{l, j}_i - \psi^{l, j}_j, w_i)
= 0, \quad \forall w_i \in \hat{W}^h_i,  \rv{\quad i=1,2,f},
\end{equation}
where
\[
W^h_i = \lbrace
w \in H^1(\omega_l): w = \delta^j_i \text{ on } \partial \omega_l
\rbrace, \quad
\hat{W}^h_i = \lbrace
w \in H^1(\omega_l): w = 0 \text{ on }\partial \omega_l
\rbrace,
\]
and $\delta^j_i$ is the piecewise constant function (delta function) for $j = 1,..,N_v^{\omega_l}$ ($N_v^{\omega_l}$ is the number of nodes on the computation mesh for $\omega_l$), $i$ is the index of continuum \rv{($i = 1,2,f$)}. Therefore, we solve $L_p^{\omega_l} = \rv{3N_v^{\omega_l}}$ local problems.

We define a snapshot space for pressures in multicontinuum media as
follows.
\begin{equation}
W_{snap}(\omega_l) = \text{span} \lbrace 
\psi^{l, j}, \, l = 1,...,N^H_v, \, j = 1,...,L_p^{\omega_l}
\rbrace.
\end{equation}
Next, we solve \rv{the} following local spectral problem on the snapshot space:
\begin{equation}
\tilde{A}_p \tilde{\phi}^l = \lambda_p \tilde{S}_p \tilde{\phi}^l, 
\end{equation}
where  $\hat{\phi}^l = (R^p_{snap})^T \tilde{\phi}^l$ and 
\[
\tilde{A}_p = R^p_{snap} {A}_p (R^p_{snap})^T, \quad 
\tilde{S}_p = R^p_{snap} {S}_p (R^p_{snap})^T, \quad 
R^p_{snap} = (\psi^{l, 1},...,\psi^{l, L_p^{\omega_l}})^T
\]
Here for matrices in triple continuum case, we have
\[
S_p = 
\begin{pmatrix}
  S_1 & 0 & 0 \\
  0 & S_2 & 0 \\
  0 & 0 & S_f 
\end{pmatrix}, \quad 
A_p = 
\begin{pmatrix}
A_1+Q_{12} + Q_{1f} & -Q_{12} & -Q_{1f} \\
-Q_{12} & A_2+Q_{12} + Q_{2f} & -Q_{2f} \\
-Q_{1f} & -Q_{2f} & A_f+Q_{1f} + Q_{2f}
\end{pmatrix}
\]
where 
\[
A_i = [a_{i,mn}], \quad 
a_{i, mn} = \int_{\omega_l^i} k_i \grad \phi^i_m \cdot \grad \phi^i_n dx, \quad 
S_i=[s_{i,mn}], \quad 
s_{i,mn}=\int_{\omega_l^i} k_i \phi^i_m \phi^i_n dx. 
\]
We choose an eigenvector $\hat{\phi}_j$ ($j = 1,..,M^{l,p}$) corresponding to the first smallest $M^{l,p}$ eigenvalues and multiply to the linear partition of unity functions $\chi^l$ for obtaining conforming basis functions
\[
W_{ms} = \text{span} \lbrace 
\phi^{l,j}, \, l = 1,...,N^H_v, \, j = 1,...,M^{l,p}
\rbrace,
\]
where $\phi^{l,j} = \chi^l \hat{\phi}^{l,j}$.

\subsection{Multiscale basis functions for displacements}
We construct the multiscale basis functions by solution following problem in local domain $\omega_l$:
find $\Psi^{l, j} \in V^h$ such that
\begin{equation}
a(\Psi^{l, j}, \bv) = 0,
\quad \forall \bv \in \hat{V}^h,
\end{equation}
where
\[
V^h = \lbrace
\bv \in H^1(\omega_l): \bv = \bar{\delta}^j_i \text{ on } \partial \omega_l
\rbrace, \quad
\hat{V}^h = \lbrace
\bv \in H^1(\omega_l): \bv = 0 \text{ on } \partial \omega_l
\rbrace.
\]
and $\bar{\delta}^j_i$ is the vector for each component for $d$-dimensional problem ($d = 2,3$) i.e. $\bar{\delta}^j_i = ({\delta}^j_i, 0, 0)$ or $\bar{\delta}^j_i = (0, {\delta}^j_i, 0)$ or $\bar{\delta}^j_i = (0, 0, {\delta}^j_i)$ for $d=3$. We solve $L_u^{\omega_l} = d \cdot N_v^{\omega_l}$ local problems.

We define snapshot space for pressures in multicontinuum media as follows
\begin{equation}
V_{snap}(\omega_l) = \text{span} \lbrace 
\Psi^{l, j}, \, l = 1,...,N^H_v, \, j = 1,...,L_u^{\omega_l}
\rbrace.
\end{equation}
For \rv{the} construction of multiscale basis, we solve \rv{the} following local spectral problem on the snapshot space:
\begin{equation}
\tilde{A}_u \tilde{\Phi} = \lambda_u \tilde{S}_u \tilde{\Phi}, 
\end{equation}
where $\hat{\Phi}^l = (R^u_{snap})^T \tilde{\Phi}^l$, 
\[
\tilde{A}_u = R^u_{snap} {A}_u (R^u_{snap})^T, \quad 
\tilde{S}_u = R^u_{snap} {S}_u (R^u_{snap})^T, \quad 
R^u_{snap} = (\Psi^{l, 1},...,\Psi^{l, L_u^{\omega_l}})^T
\]
and 
\[
A_u = [a_{u,mn}], \quad  
a_{u,mn} = \int_{\omega_l} \bsig(\Phi_m) \cdot \bvar(\Phi_n) dx,
\quad 
S_u=[s_{i,mn}], \quad 
s_{i,mn}=\int_{\omega_l} (\lambda + 2\mu) \Phi^i_m \Phi^i_n dx. 
\]
We choose \rv{eigenvectors} $\hat{\Phi}_j$, $j = 1,..,M^{l,u}$
corresponding to the first smallest $M^{l,u}$ eigenvalues and multiply
\rv{by}
the linear partition of unity functions \rv{to obtain the} conforming basis functions:
\[
V_{ms} = \text{span} \lbrace 
\Phi^{l,j}, \, l = 1,...,N^H_v, \, j = 1,...,M^{l,u}
\rbrace,
\]
where $\Phi^{l,j} = \chi^l \hat{\Phi}^{l,j}$.

\subsection{Coarse grid system}. 
Using \rv{the above} constructed multiscale basis functions for
pressures and displacements, we define \rv{the} projection matrix:
\begin{equation}
R = \begin{pmatrix}
R_p & 0 \\
0 & R_u
\end{pmatrix},
\end{equation}
where 
\[
R_u = (\Phi^{1,1}, ... , \Phi^{1,M^{1,u}} , ..., \Phi^{N^H_v,1}, ...,\Phi^{N^H_v, M^{N^H_v,u}})^T,
\]\[
R_p = (\phi^{1,1}, ... , \phi^{1,M^{1,p}} , ..., \phi^{N^H_v,1}, ...,\phi^{N^H_v, M^{N^H_v,p}})^T.
\]
Then we obtain the following reduced order model:
\begin{equation}
\begin{split}
\zeta^{(\alpha_i)}_{\tau} M^H_i p^{H,n}_i + 
& \zeta^{(\beta_i)}_{\tau} D^H_i \bu^{H,n} + 
A^H_i p^{H,n}_i + 
\sum_{j!=i} Q^H_{ij}(p^{H,n}_i - p^{H,n}_j) \\
&= 
\zeta^{(\alpha_i)}_{\tau} M^H_i p^{H,n-1}_i -
\zeta^{(\alpha_i)}_{\tau} 
\sum_{j=2}^{n} \zeta^{(\alpha_i)}_{j-1} M^H_i (p^{H, n-j+1}_i - p^{H, n-j}_i) \\
& + 
\zeta^{(\beta_i)}_{\tau} D^H_i \bu^{H, n-1} -
\zeta^{(\beta_i)}_{\tau}
\sum_{j=2}^{n} \zeta^{(\beta_i)}_{j-1} D^H_i (\bu^{H, n-j+1} - \bu^{H, n-j}), 
\\
&\sum_{j} (D^H_j)^T p^{H,n}_j + A^H_u \bu^{H,n} = 0,
\end{split}
\end{equation}
where 
\[
M^H_i = R M_i R^T, \quad 
A^H_i = R A_i R^T, \quad 
Q^H_{ij} = R Q_{ij} R^T, \quad 
D^H_i = R D_i R^T, \quad 
A^H_u = R A_u R^T.
\]
After obtaining \rv{the coarse-scale solutions}, we reconstruct \rv{the} fine-scale solutions:
\[
p^{ms,n}_i = R^T p^{H,n}_i,  \quad
u^{ms,n} = R^T u^{H,n}.
\]
\rv{We remark that in our method presented above we store and use only the
information on the coarse-grid solutions at the previous time step.}

\section{Numerical results}
\begin{figure}[h!]
\centering
\includegraphics[width=0.30\textwidth]{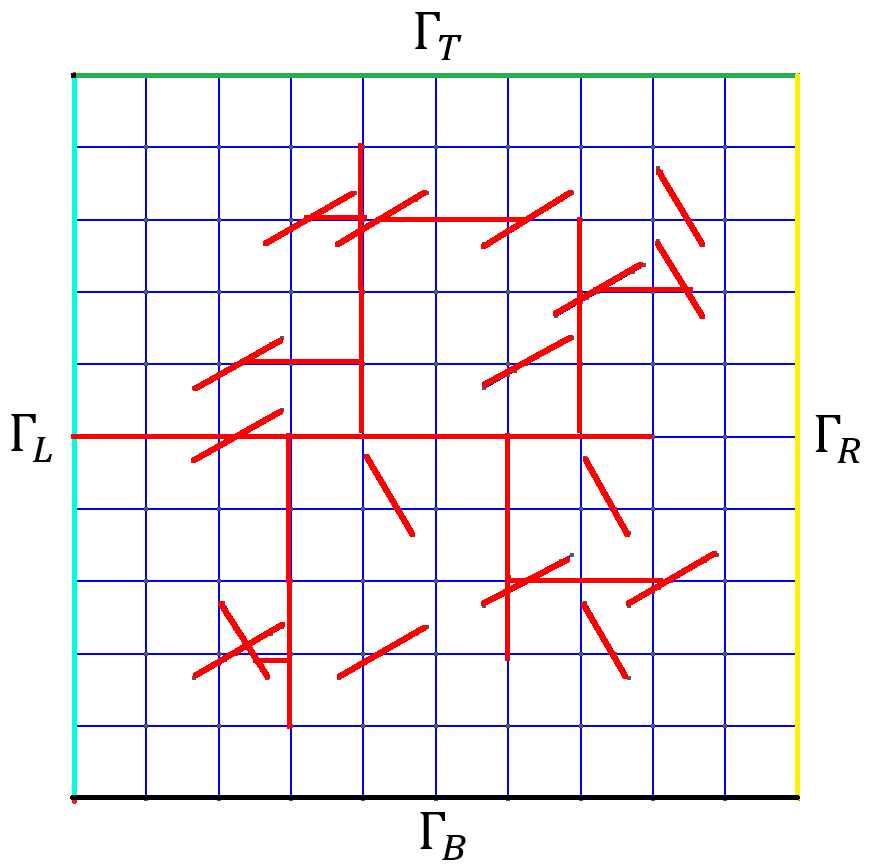}
\includegraphics[width=0.30\textwidth]{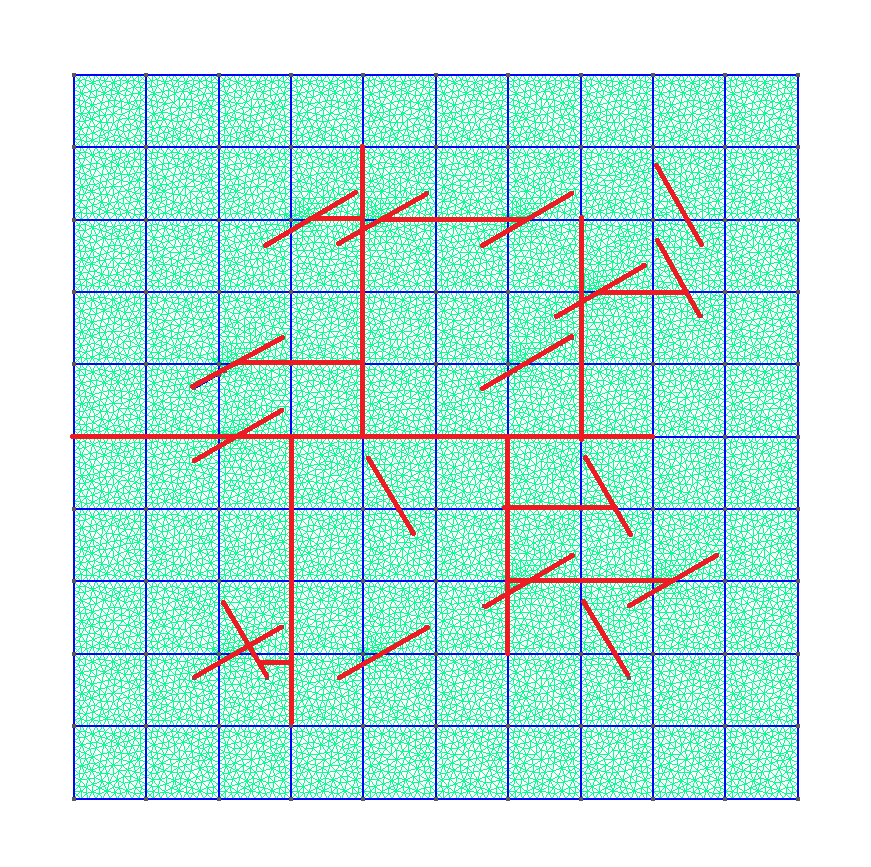}
\caption{Computation domain and grids.  Coarse grid (blue color), fine grid (green), and fractures (red).}
\label{mesh2frac}
\end{figure}

\begin{figure}[h!]
\centering
\includegraphics[width=0.30\linewidth]{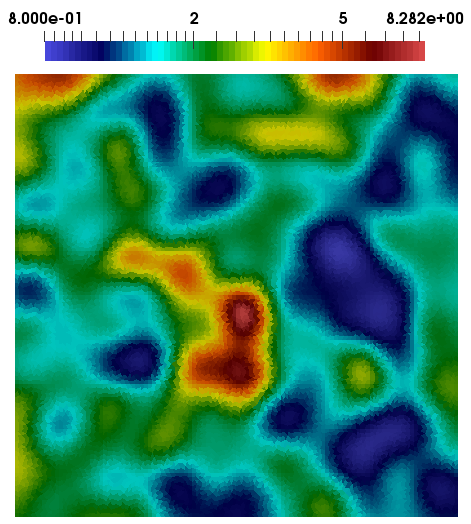}
\includegraphics[width=0.30\linewidth]{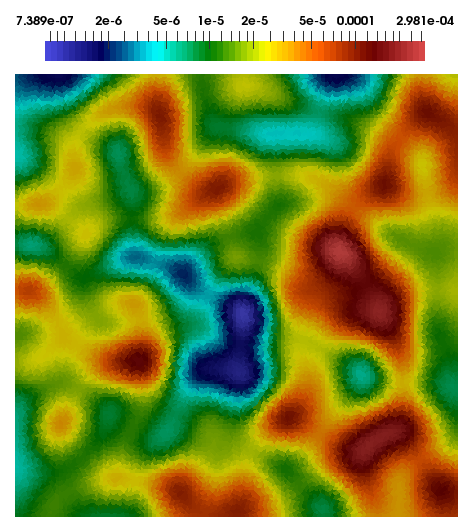}
\includegraphics[width=0.30\linewidth]{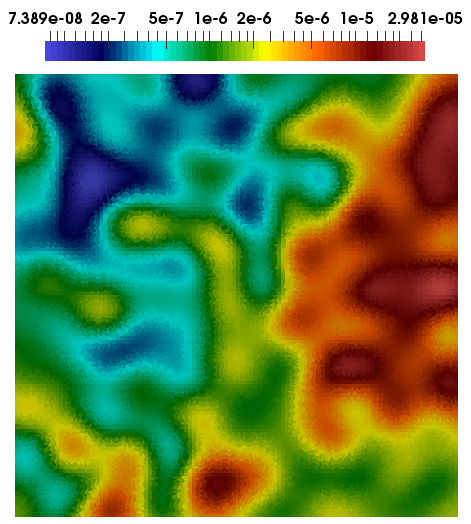}
\caption{Elasticity coefficient $E$(left) and heterogeneous permeability $k_1$(center) and $k_2$(right).}
\label{hetero}
\end{figure}

In this section, we present the numerical results of the
poroelasticity problems in heterogeneous and fractured media with
fractional derivatives. The coarse grid is uniform with rectangular
cells. In Figure \ref{mesh2frac}, we show computational coarse and
fine grids. The fine grid contains 25846 cells and 12944 vertices, and
the coarse grid contains 121 vertices and 100 cells. We consider the
time-fractional diffusion equation for poroelasticity problem in
$\Omega = (0,50)^2$ for two cases such as the poroelasticity in
fractured media and multicontinuum media.
For coefficients representing matrix and fracture properties,
we set $\gamma_1 = 0.1, \gamma_2 = 0.1, \gamma_f = 0, k_f = 1.0, M_1 =
M_2 = 10, M_f = 10^3, \nu = 0.3.$
The calculation is performed by $T_{max} = 86400$ with times step
$\tau = 8640$ and $\eta_{12} = 5*k_2$.
Heterogeneous coefficients for elasticity modulus and heterogeneous
permeability for \rv{the first and second continua} are presented in
Fig. \ref{hetero}. \rv{A numerical solution is presented with} the
following boundary conditions
$u_x = 0, \bsig_y = 0, x \in \Gamma_L \cup \Gamma_R$, $u_y = 0,
\bsig_x =0, x \in \Gamma_T \cup \Gamma_B$ for displacement, and the initial condition $p^0 = 1$ for pressure.

To compare the results, we use the fine grid solution as a reference
solution and calculate relative $L_2$ norm and $H_1$ semi-norm of
errors between \rv{the multiscale and reference solutions}
\[
e_{L^2}^{p_i}  =  \rv{ \left(
\frac{(p_i - p_i^{ms}, p_i-p_i^{ms})}{(p_i, p_i)} \right)^{1/2}, }
\quad  
e^u_{L^2} = \rv{ \left( 
\frac{(\bu - \bu^{ms}, \bu-\bu^{ms})}{(\bu,\bu)} 
\right)^{1/2}, }
\]
\[
e_{H_1}^{p_i}  =   \left(
\frac{b(p_i - p_i^{ms}, p_i-p_i^{ms})}{b(p_i, p_i)} \right)^{1/2}, \quad  
e^u_{H_1} =  \left( 
\frac{a(\bu - \bu^{ms}, \bu-\bu^{ms})}{a(\bu,\bu)} 
\right)^{1/2}, 
\]
where $i$ is the index for the continuum $(i = 1,2),$ and
$y^{ms} = (p_1^{ms}, p_2^{ms}, \bu^{ms})$ denotes the multiscale solution
using the GMsFEM and $y = (p_1, p_2, \bu)$ the fine grid solution.


\subsection{Poroelasticity in fractured media}

We present how the introduction of the time memory effect, by means of the Caputo’s fractional time derivative in the constitutive equation, affects both the pressure and displacement in fractured media. In this subsection, we solve the poroelasticity problem with one continuum. 

\begin{figure}[h!]
\centering
\vspace{2.0 mm}
\includegraphics[width=0.25\linewidth]{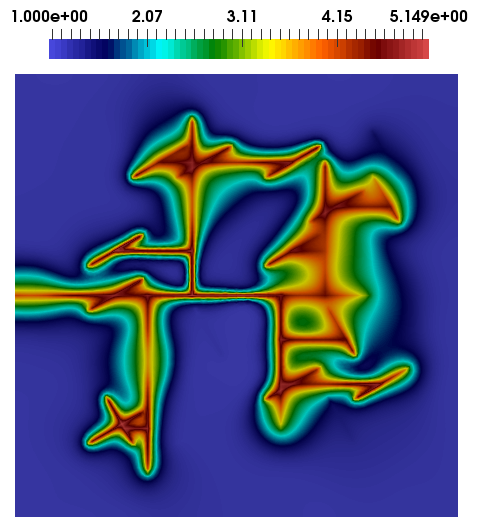}
\includegraphics[width=0.25\linewidth]{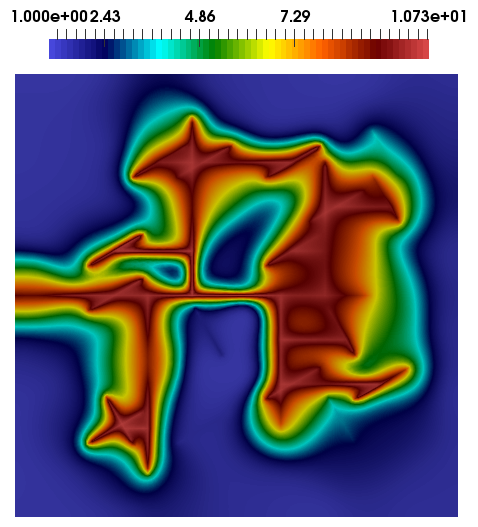}
\includegraphics[width=0.25\linewidth]{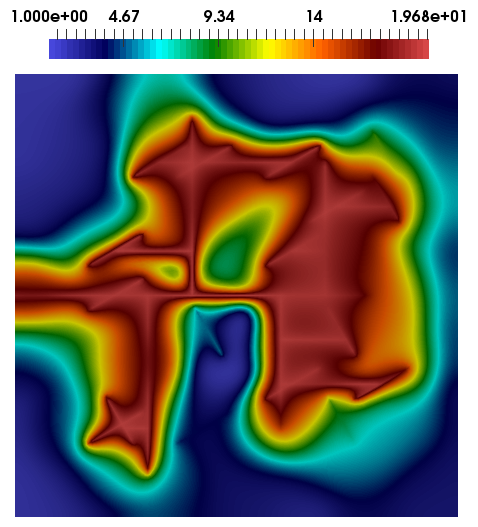}\\
\includegraphics[width=0.25\linewidth]{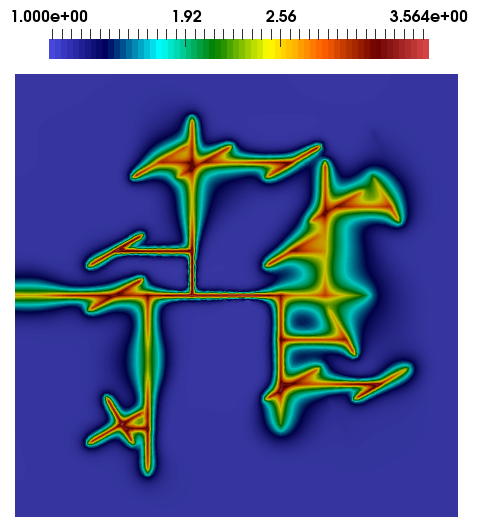}
\includegraphics[width=0.25\linewidth]{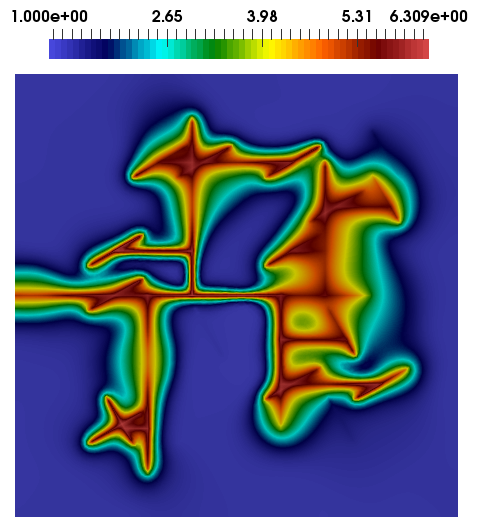}
\includegraphics[width=0.25\linewidth]{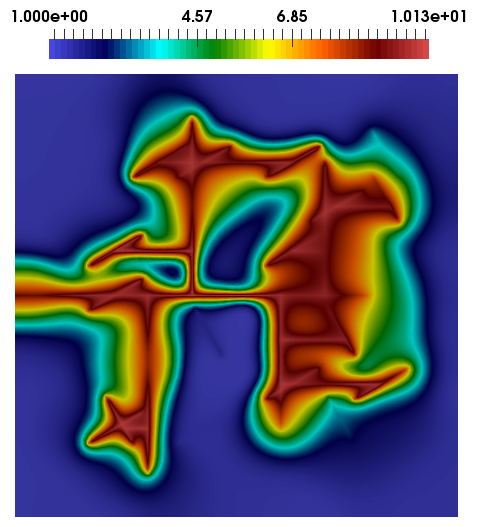}\\
\includegraphics[width=0.25\linewidth]{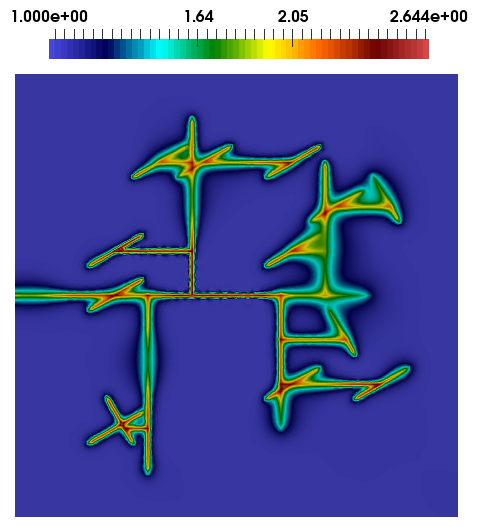}
\includegraphics[width=0.25\linewidth]{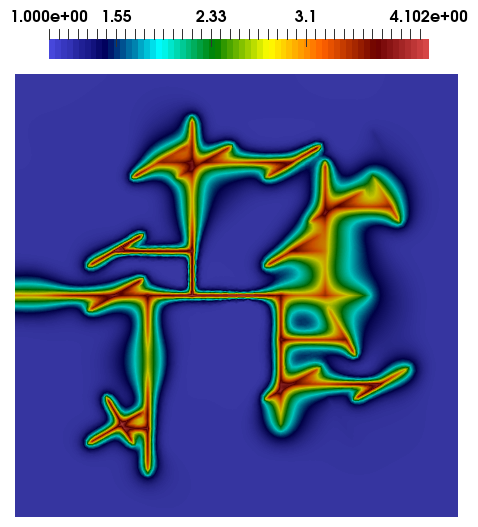}
\includegraphics[width=0.25\linewidth]{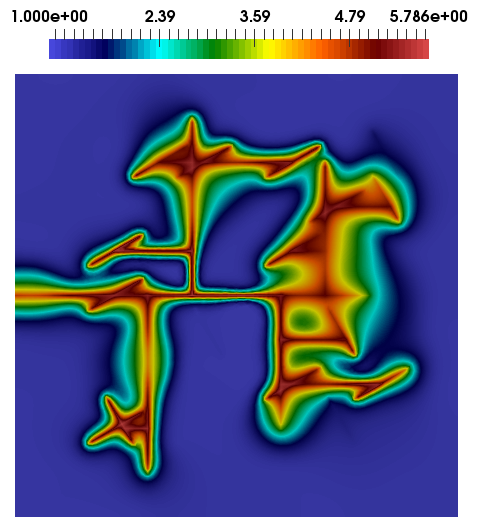}
\caption{Distribution of pressure for exact solution for $t=0$, $25920$ and $t=86400$ (from left to right). 
First row: fractional order derivative $\alpha = 1.0$.  
Second row: fractional order derivative $\alpha = 0.9$. 
Third row: fractional order derivative $\alpha = 0.8$.  }
\label{upexact}
\end{figure}

In Figure $\ref{upexact}$ we present the numerical solution
distribution of pressure at different time steps on a fine mesh with
other fractional order derivative. Relative $L_2$ and energy $H_1$
errors are presented for different number of multiscale basis
functions in Tables $\ref{table1}-\ref{table3}$. We present the error
comparison between the fine-scale and multiscale solutions with different numbers of multiscale basis functions. We observe that the error decreases when we increase the number of multiscale basis functions for each fractional order derivative. The relative error reduces from $21\%$ to $0.8\%$ for displacement and $12\%$ to $0.4\%$ for pressure with fractional order derivative $\alpha = 1.0$. To obtain a good solution, we need to take twelve basis functions in each fractional order derivative.
 
\begin{figure}[h!]
\centering
\includegraphics[width=0.40\textwidth]{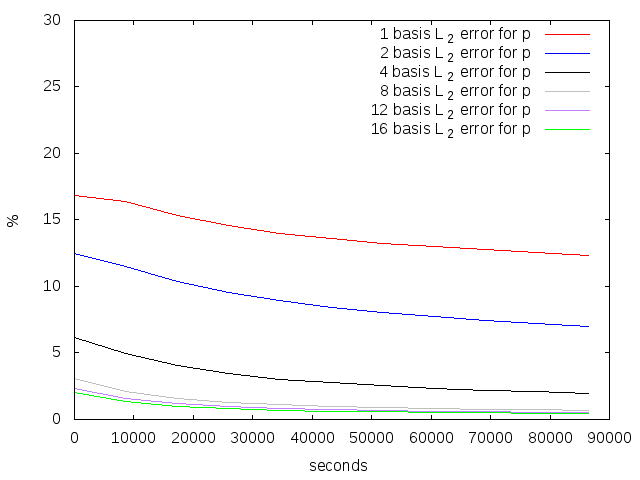}
\includegraphics[width=0.40\textwidth]{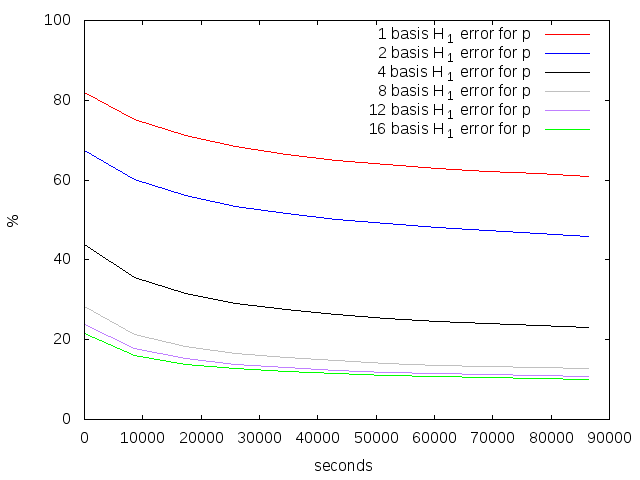}\\
\includegraphics[width=0.40\textwidth]{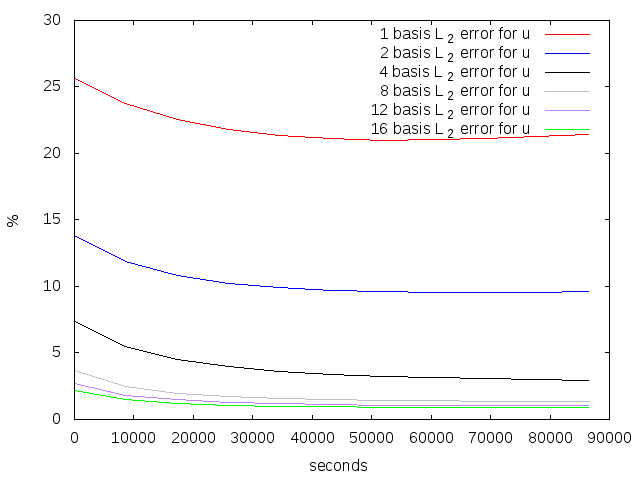}
\includegraphics[width=0.40\textwidth]{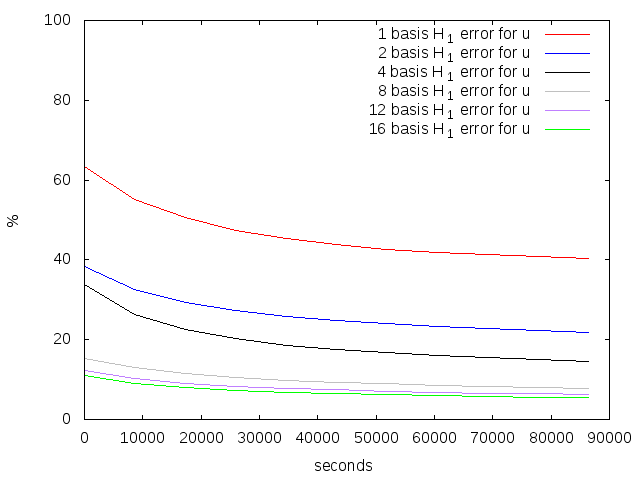}
\caption{Poroelasticity in fractured media: Relative $L_2$ error (left) and $H_1$ (right) errors vs time for different number multiscale basis functions for pressure (first row) and displacements (second row) with fractional order derivative $\alpha = 1.0$. 
 }
\label{graphic-1}
\end{figure}

\begin{figure}[h!]
\centering
\includegraphics[width=0.40\textwidth]{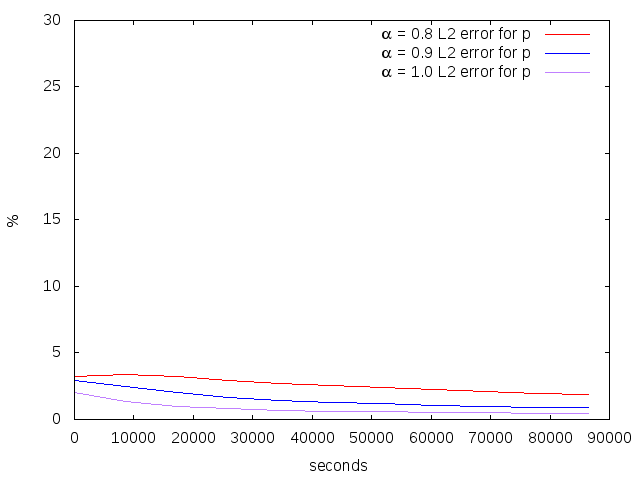}
\includegraphics[width=0.40\textwidth]{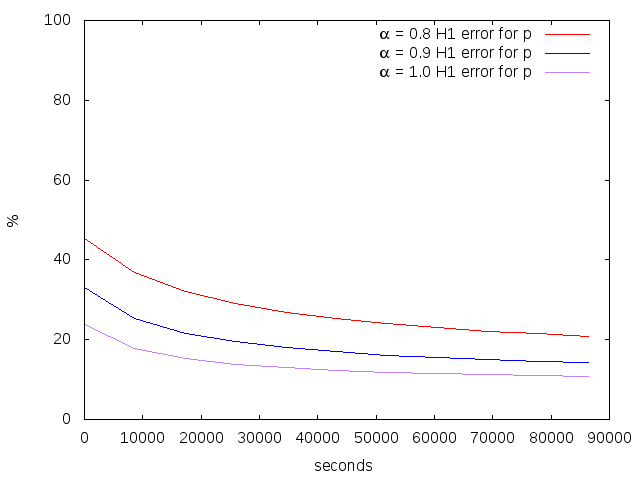}\\
\includegraphics[width=0.40\textwidth]{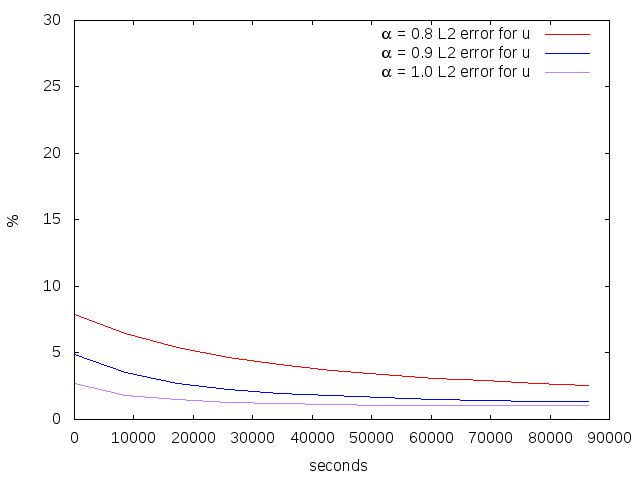}
\includegraphics[width=0.40\textwidth]{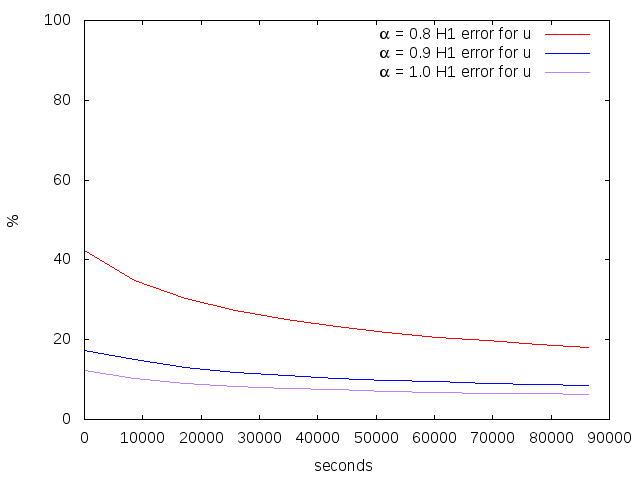}
\caption{Poroelasticity in fractured media: Relative $L_2$ error (left) and $H_1$ (right) errors vs time for different fractional order derivative($\alpha = 0.8,0.9,1.0$) for pressure (first row) and displacements (second row) with multiscale basis function 12. 
 }
\label{graphic-2}
\end{figure}

Next, the relative $L_2$ error dynamics in $\%$ for different number of multiscale basis functions with fractional order derivative $\alpha = 1.0$ are shown in Figure $\ref{graphic-1}$. We observe that the errors reduce by time. In Figure $\ref{graphic-2}$ we present relative $L_2$ error dynamics with different fractional order derivative for twelve multiscale basis functions. The behavior of the figures is similar to the previous figures. Therefore, we can assume that the method provides a good solution.

\begin{figure}[h!]
\centering
\vspace{2.0 mm}
\includegraphics[width=0.25\linewidth]{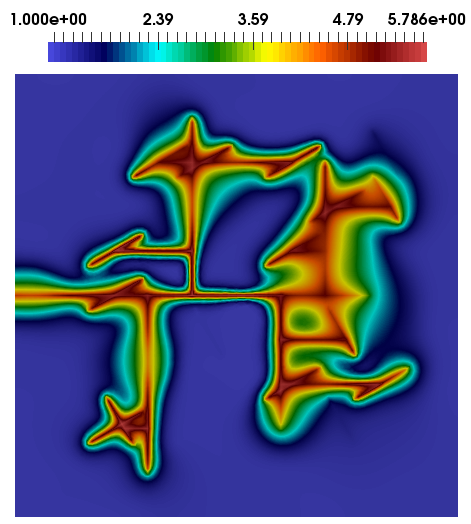}
\includegraphics[width=0.25\linewidth]{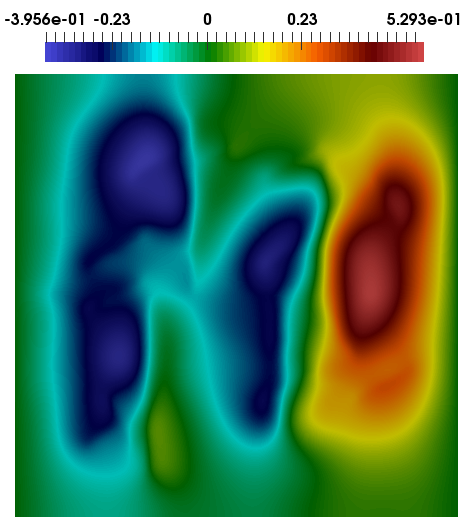}  
\includegraphics[width=0.25\linewidth]{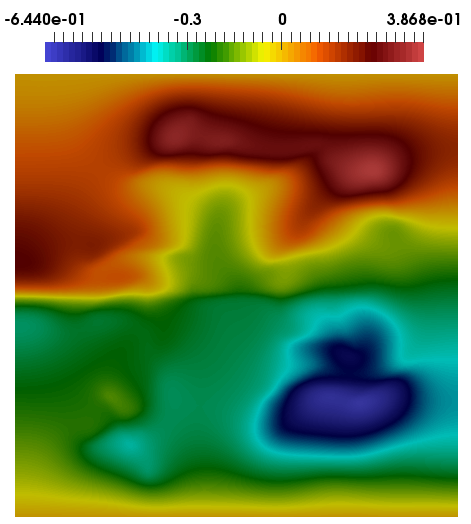}\\
\includegraphics[width=0.25\linewidth]{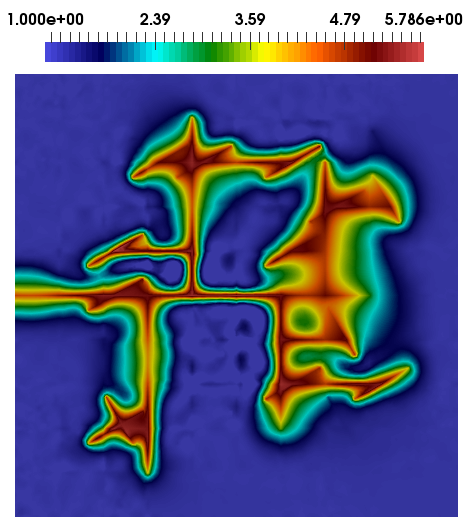} 
\includegraphics[width=0.25\linewidth]{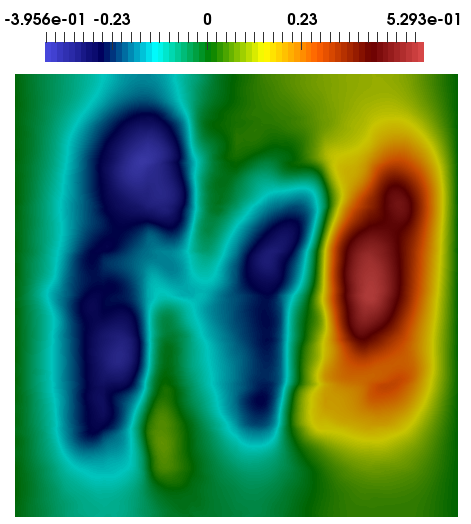}
\includegraphics[width=0.25\linewidth]{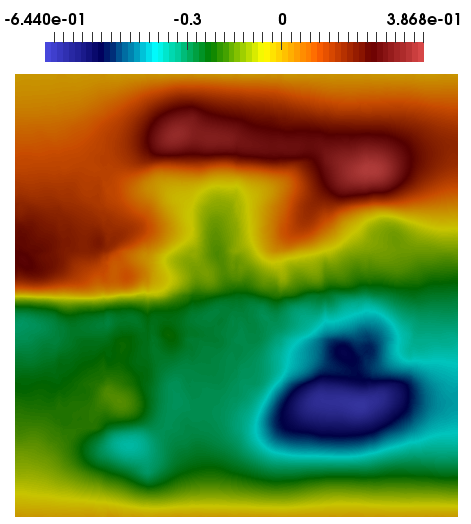} 
\caption{Poroelasticity in fractured media: Distribution of pressure and displacement along $X$ and $Y$ (from left to right) at final time for fractional order derivative $\alpha = 0.8$.
First row: exact solution. Second row: multiscale solution(12 multiscale basis functions).}
\label{pic1}
\end{figure}

\begin{figure}[h!]
\centering
\vspace{2.0 mm}
\includegraphics[width=0.25\linewidth]{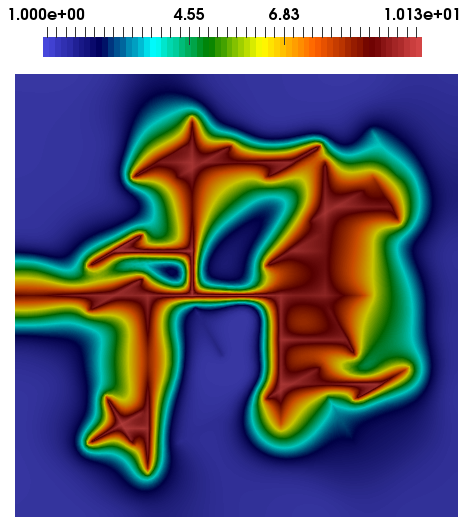}
\includegraphics[width=0.25\linewidth]{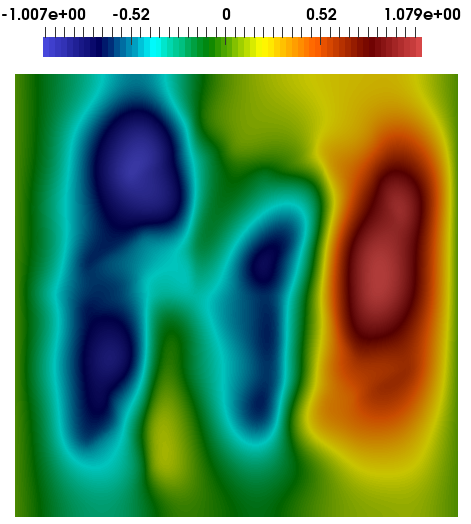}  
\includegraphics[width=0.25\linewidth]{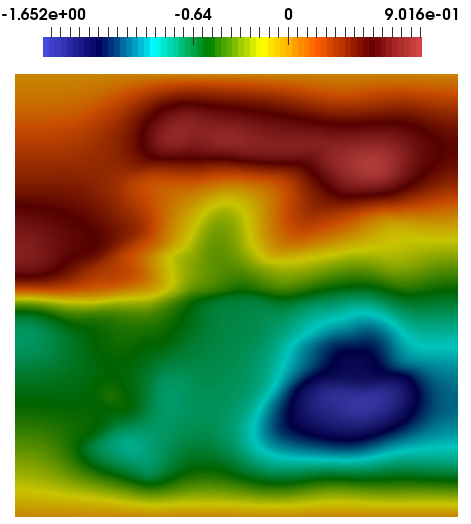}\\
\includegraphics[width=0.25\linewidth]{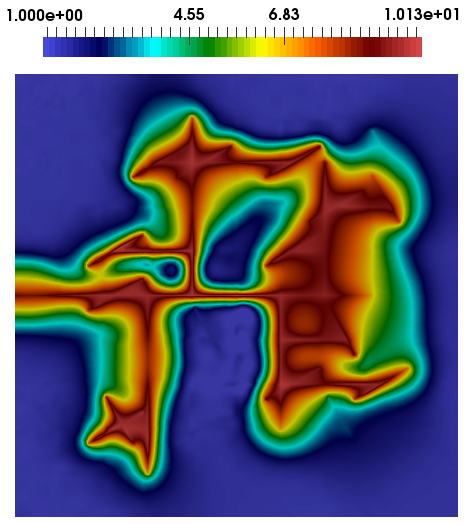} 
\includegraphics[width=0.25\linewidth]{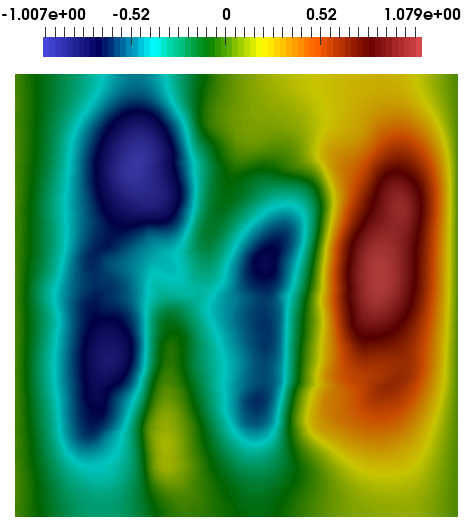}
\includegraphics[width=0.25\linewidth]{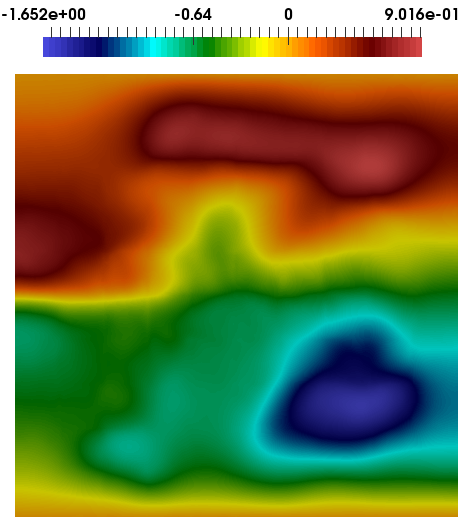} 
\caption{Poroelasticity in fractured media: Distribution of pressure and displacement along $X$ and $Y$ (from left to right) at final time for fractional order derivative $\alpha = 0.9$.
First row: exact solution. Second row: multiscale solution(12 multiscale basis functions).}
\label{pic2}
\end{figure}

\begin{figure}[h!]
\centering
\vspace{2.0 mm}
\includegraphics[width=0.25\linewidth]{1exactp10}
\includegraphics[width=0.25\linewidth]{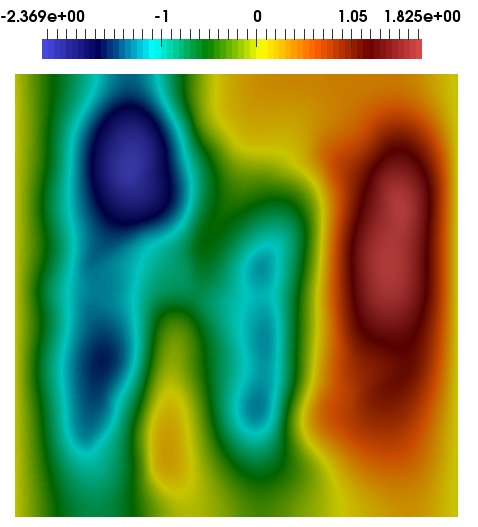}  
\includegraphics[width=0.25\linewidth]{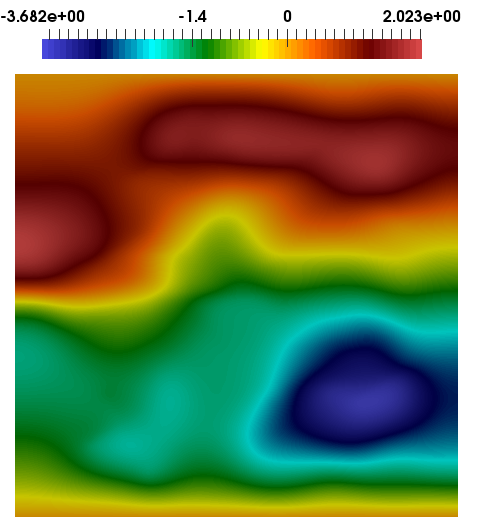}\\
\includegraphics[width=0.25\linewidth]{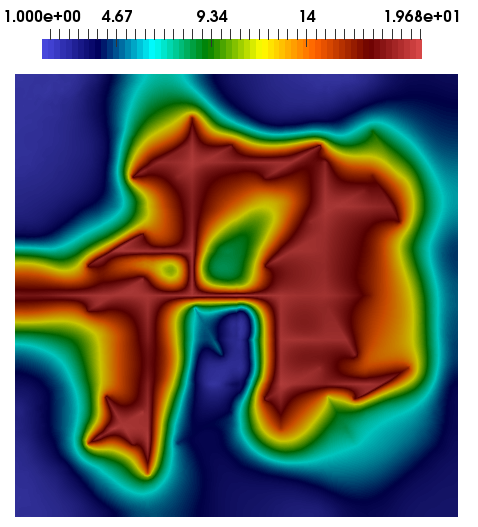} 
\includegraphics[width=0.25\linewidth]{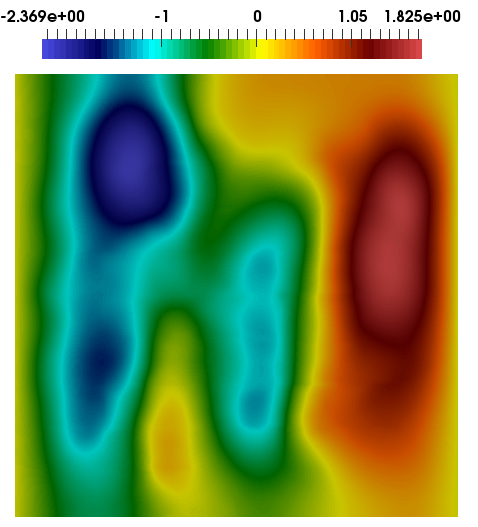}
\includegraphics[width=0.25\linewidth]{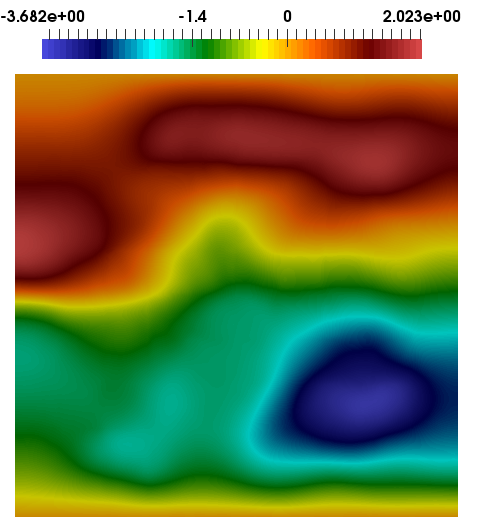} 
\caption{Poroelasticity in fractured media: Distribution of pressure and displacement along $X$ and $Y$ (from left to right) at final time for fractional order derivative $\alpha = 1.0$.
First row: exact solution. Second row: multiscale solution(12 multiscale basis functions).}
\label{pic3}
\end{figure}

Figures $\ref{pic1}-\ref{pic3}$ \rv{show} the distribution of pressure
and displacement along $X$ and $Y$ directions at final time for
different fractional order derivatives. In the first row, we show fine
scale and multiscale solutions with twelve multiscale basis functions
for the GMsFEM is presented in the second row. We observe good results of the presented method for solving poroelasticity problems for different fractional order derivatives.

\begin{table}[h!]
\begin{center}
\begin{tabular}{ | c | c | c | c | c | c | }
\hline
$M$ & $DOF_H$ & $e_{L_2}^u$ (\%) & $e_{H_1}^u$  (\%) & $e_{L_2}^p$  (\%) & $e_{H_1}^p$  (\%) \\ \hline
1 & 363 & 26.324 & 62.768 & 18.280 & 80.061\\
2 & 726 & 14.436 & 40.353 & 13.540 & 64.607\\
4 & 1452 & 7.355 & 30.629 & 6.734 & 40.175\\
8 & 2904 & 3.641 & 21.532 & 3.149 & 25.884\\
12 & 4356 & 2.580 & 18.158 & 2.261 & 20.837\\
16 & 5808 & 2.112 & 16.312 & 1.896 & 18.730\\
\hline
\end{tabular}
\end{center}
\caption{Relative errors for displacement and pressure with different numbers of multiscale basis functions. Poroelasticity in fractured media with fractional order derivative $\alpha = 0.8$}
\label{table1}
\end{table}

\begin{table}[h!]
\begin{center}
\begin{tabular}{ | c | c | c | c | c | c | }
\hline
$M$ & $DOF_H$ & $e_{L_2}^u$ (\%) & $e_{H_1}^u$  (\%) & $e_{L_2}^p$  (\%) & $e_{H_1}^p$  (\%) \\ \hline
1 & 363 & 21.975 & 48.612 & 14.909 & 69.075\\
2 & 726 & 10.425 & 28.102 & 9.898 & 54.202\\
4 & 1452 & 4.112 & 20.956 & 3.651 & 29.799\\
8 & 2904 & 1.803 & 10.888 & 1.406 & 17.002\\
12 & 4356 & 1.333 & 8.603 & 1.027 & 14.209\\
16 & 5808 & 1.115 & 7.463 & 0.872 & 13.023\\
\hline
\end{tabular}
\end{center}
\caption{Relative errors for displacement and pressure with different numbers of multiscale basis functions. Poroelasticity in fractured media with fractional order derivative $\alpha = 0.9$}
\label{table2}
\end{table}

\begin{table}[h!]
\begin{center}
\begin{tabular}{ | c | c | c | c | c | c | }
\hline
$M$ & $DOF_H$ & $e_{L_2}^u$ (\%) & $e_{H_1}^u$  (\%) & $e_{L_2}^p$  (\%) & $e_{H_1}^p$  (\%) \\ \hline
1 & 363 & 21.461 & 40.372 & 12.348 & 61.024\\
2 & 726 & 9.611 & 21.865 & 6.983 & 45.957\\
4 & 1452 & 2.959 & 14.532 & 1.979 & 22.994\\
8 & 2904 & 1.341 & 7.890 & 0.680 & 12.675\\
12 & 4356 & 1.017 & 6.256 & 0.510 & 10.780\\
16 & 5808 & 0.871 & 5.402 & 0.450 & 10.063\\
\hline
\end{tabular}
\end{center}
\caption{Relative errors for displacement and pressure with different numbers of multiscale basis functions. Poroelasticity in fractured media with fractional order derivative $\alpha = 1.0$}
\label{table3}
\end{table}

\subsection{Poroelasticity in multicontinuum media}
Next we consider a computational macroscopic model for fractional Biot's system in multicontinuum heterogeneous media. In this case, we consider two continua. In this part of the work, we investigate a time memory formalism in poroelasticity problem in multicontinuum media.

In Tables $\ref{tablemult1}$-$\ref{tablemult3}$ relative $L_2$ and
$H_1$-energy  errors are presented for different number of multiscale
basis functions. The results show that twelve multiscale basis
functions are enough to achieve good results, for example, for
fractional order derivative $\alpha = 1.0$ with $1.051 \%$ of $L_2$
error for displacement, $0.667 \%$ and $0.644 \%$ of $L_2$ errors for the first and second continuum pressures. We have similar improvements for further increment of the multiscale basis functions. 

\begin{figure}[h!]
\centering
\includegraphics[width=0.40\textwidth]{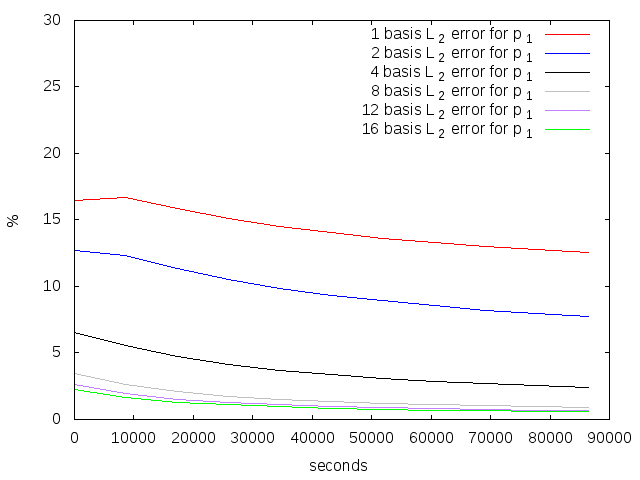}
\includegraphics[width=0.40\textwidth]{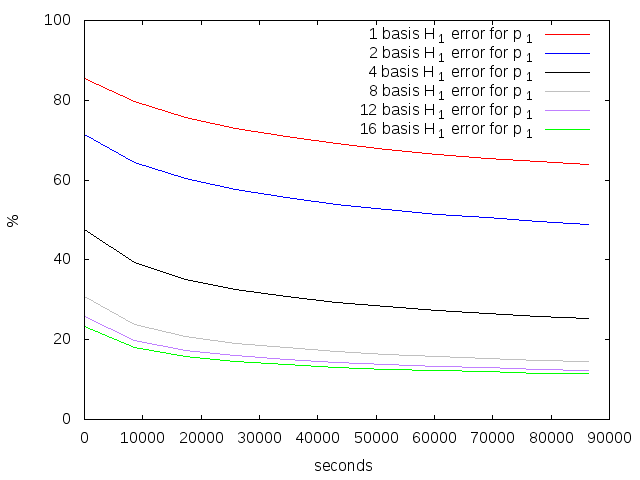}\\
\includegraphics[width=0.40\textwidth]{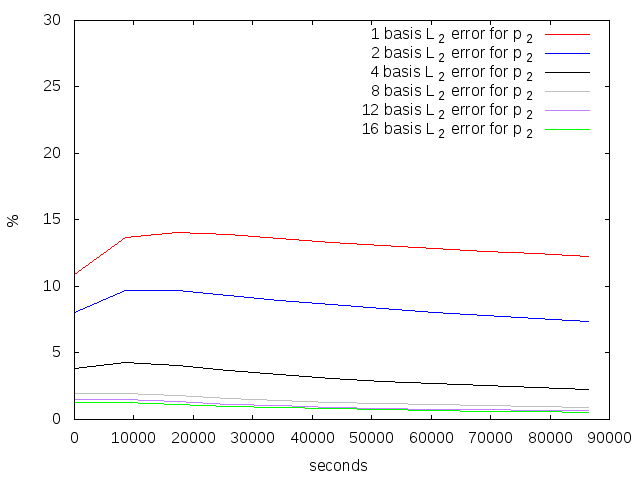}
\includegraphics[width=0.40\textwidth]{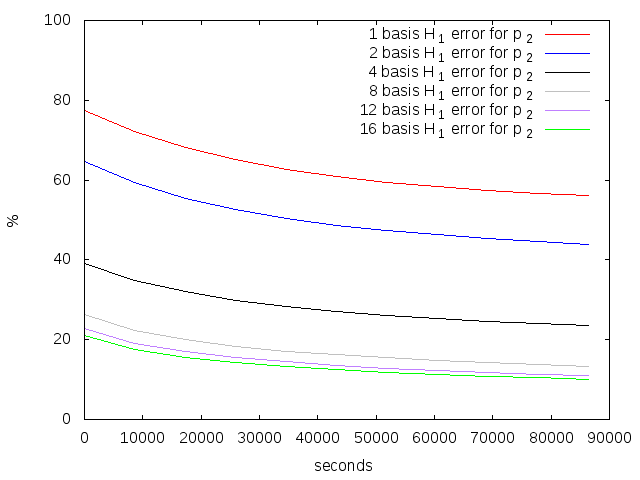}\\
\includegraphics[width=0.40\textwidth]{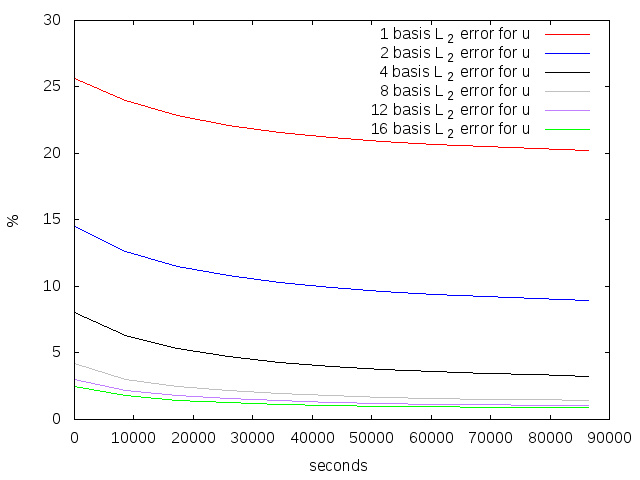}
\includegraphics[width=0.40\textwidth]{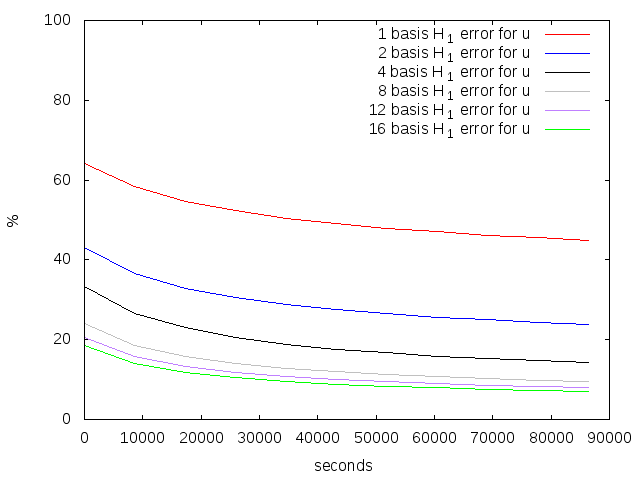}
\caption{Poroelasticity in multicontinuum media: Relative $L_2$ error (left) and $H_1$ (right) errors vs time for different number multiscale basis functions for first continuum (first row), second continuum (second row) and displacements (third row) with fractional order derivative $\alpha = 1.0$. 
 }
\label{graphic-3}
\end{figure}

\begin{figure}[h!]
\centering
\includegraphics[width=0.40\textwidth]{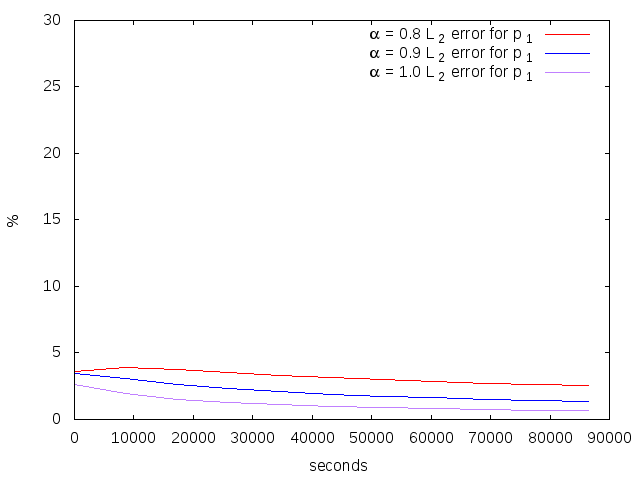}
\includegraphics[width=0.40\textwidth]{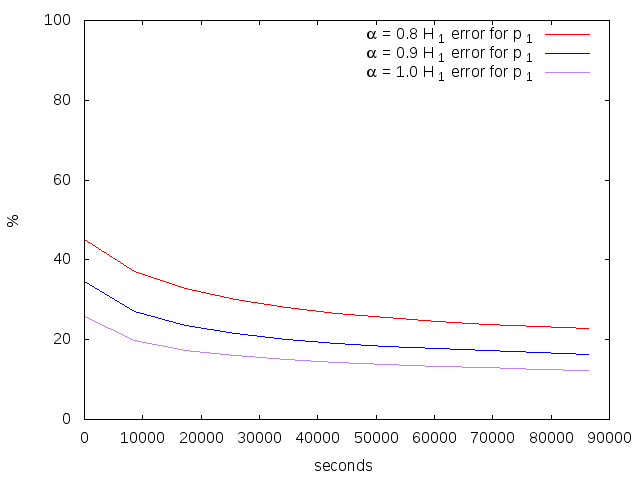}\\
\includegraphics[width=0.40\textwidth]{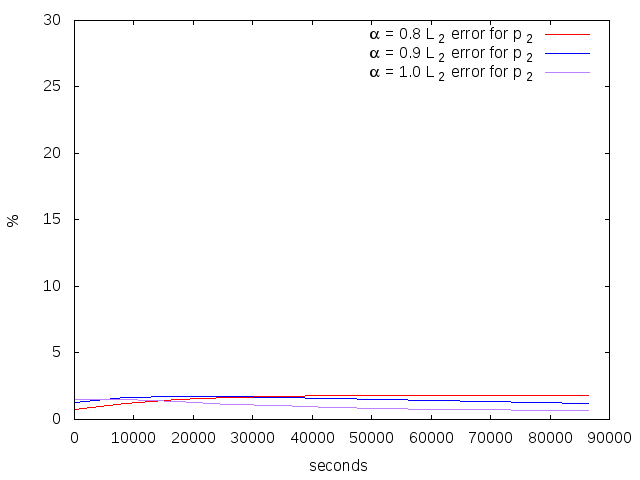}
\includegraphics[width=0.40\textwidth]{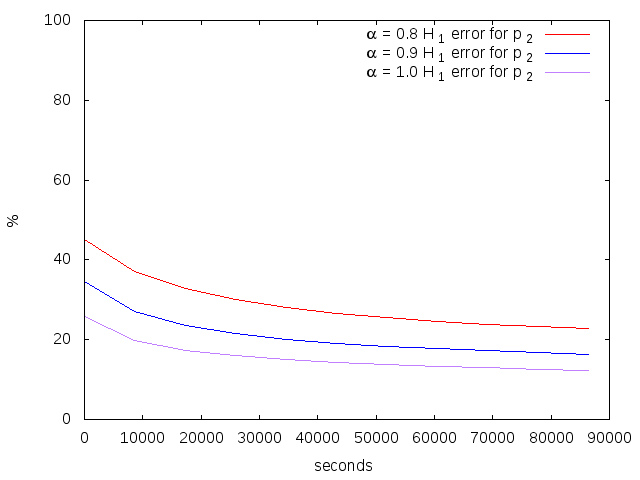}\\
\includegraphics[width=0.40\textwidth]{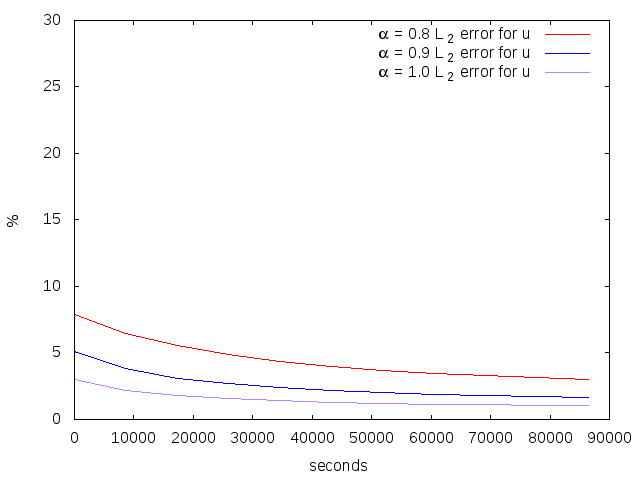}
\includegraphics[width=0.40\textwidth]{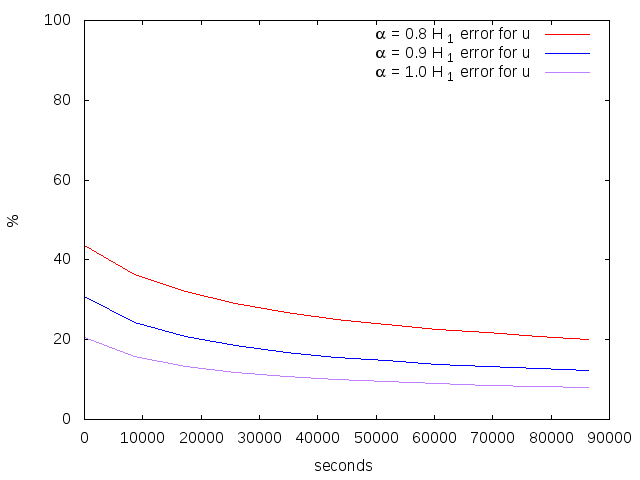}
\caption{Poroelasticity in multicontinuum media: Relative $L_2$ error (left) and $H_1$ (right) errors vs time for different fractional order derivative($\alpha = 0.8,0.9,1.0$) for first continuum (first row), second continuum (second row) and displacements (third row) with multiscale basis function 12.}
\label{graphic-4}
\end{figure}

Then the relative $L_2$ error dynamics in $\%$ for different number of multiscale basis functions with fractional order derivative $\alpha = 1.0$ are shown in Figure $\ref{graphic-3}$. We also observe that errors reduce by time for poroelasticity in multicontinuum media. In Figure $\ref{graphic-4}$ we present relative $L_2$ error dynamics with different fractional order derivative for twelve multiscale basis functions. We observe that the presented method provides good results for different fractional order derivatives.

\begin{figure}[h!]
\centering
\vspace{2.0 mm}
\includegraphics[width=0.24\linewidth]{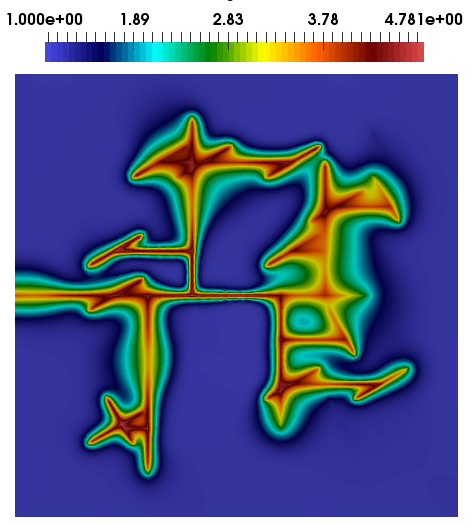}
\includegraphics[width=0.24\linewidth]{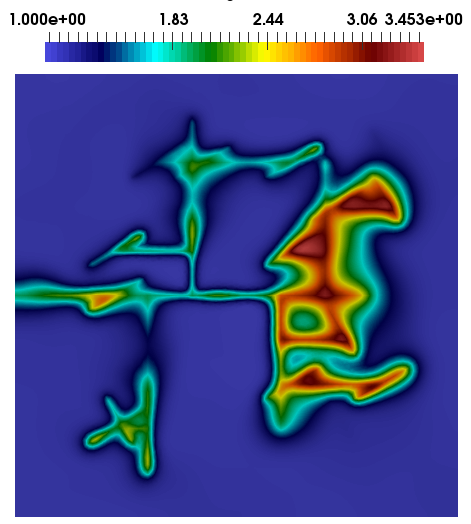}  
\includegraphics[width=0.24\linewidth]{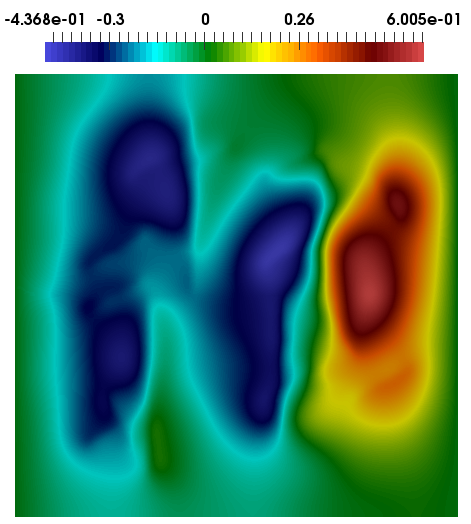}
\includegraphics[width=0.24\linewidth]{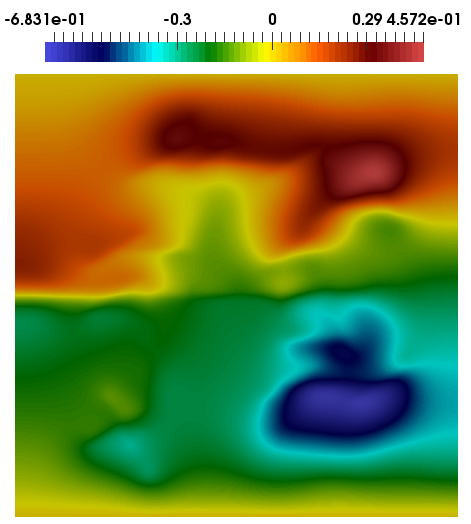} \\
\includegraphics[width=0.24\linewidth]{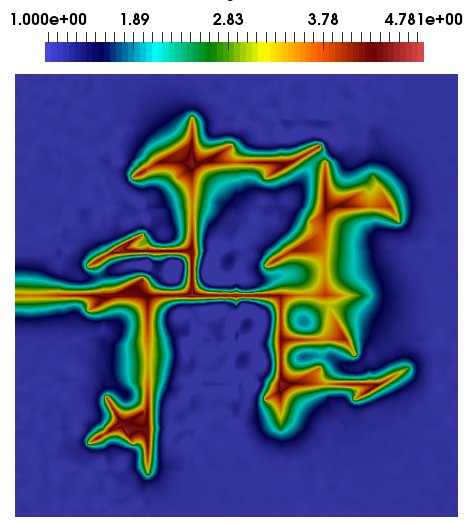}
\includegraphics[width=0.24\linewidth]{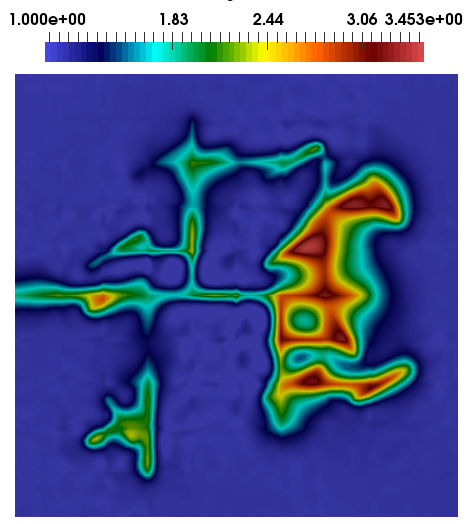}
\includegraphics[width=0.24\linewidth]{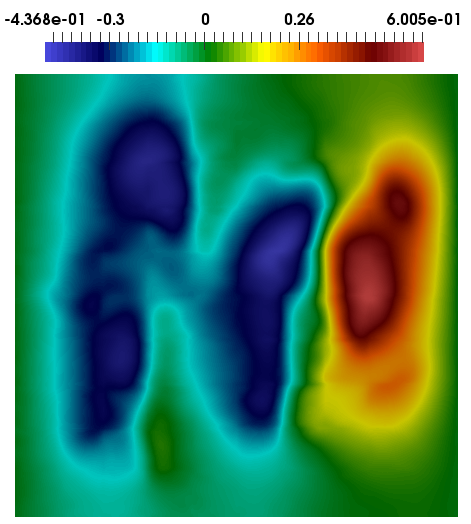}
\includegraphics[width=0.24\linewidth]{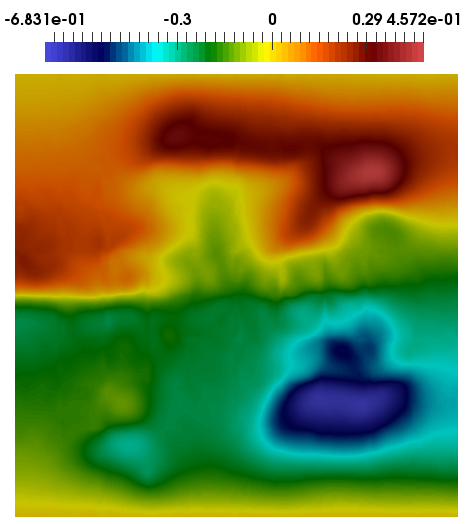} \\
\caption{Poroelasticity in multicontinuum media: Distribution of pressure for first and second continuum and displacement along $X$ and $Y$ at final time for fractional order derivative $\alpha = 0.8$ (from left to right).
First row: exact solution. Second row: multiscale solution(12 multiscale basis functions).}
\label{pic-4}
\end{figure}

\begin{figure}[h!]
\centering
\vspace{2.0 mm}
\includegraphics[width=0.24\linewidth]{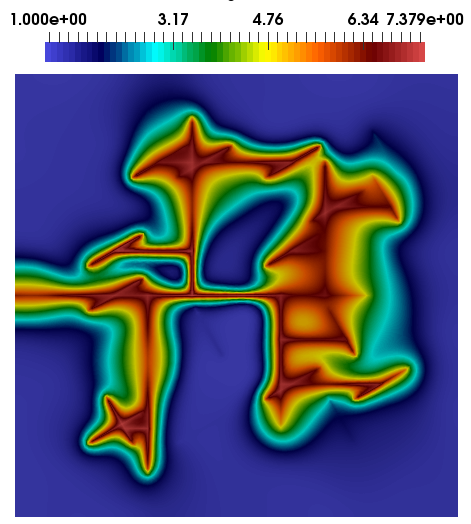}
\includegraphics[width=0.24\linewidth]{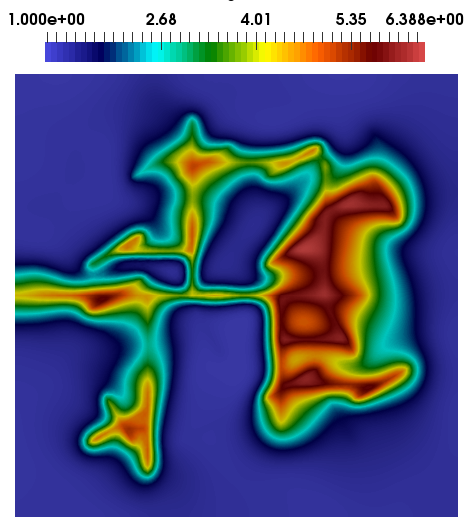}  
\includegraphics[width=0.24\linewidth]{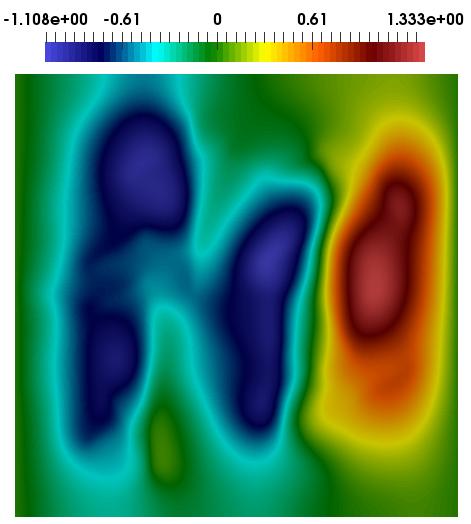}
\includegraphics[width=0.24\linewidth]{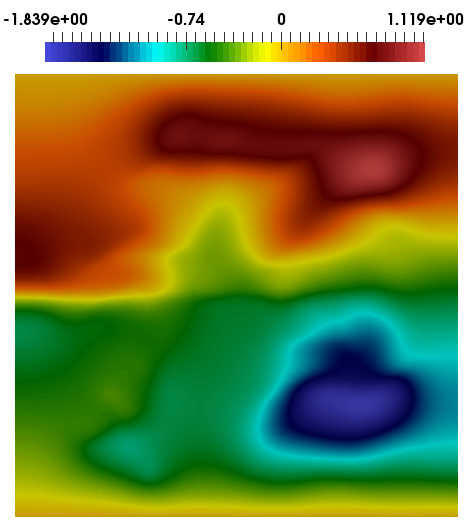} \\
\includegraphics[width=0.24\linewidth]{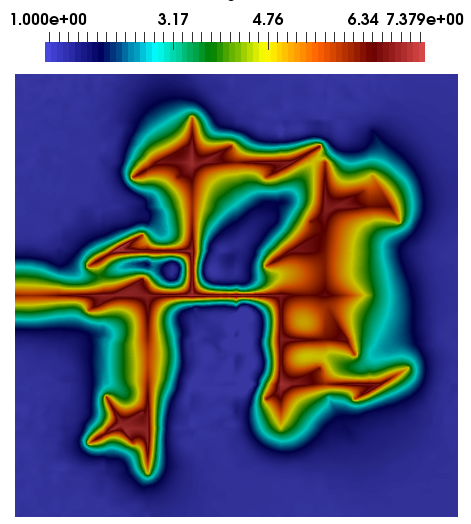}
\includegraphics[width=0.24\linewidth]{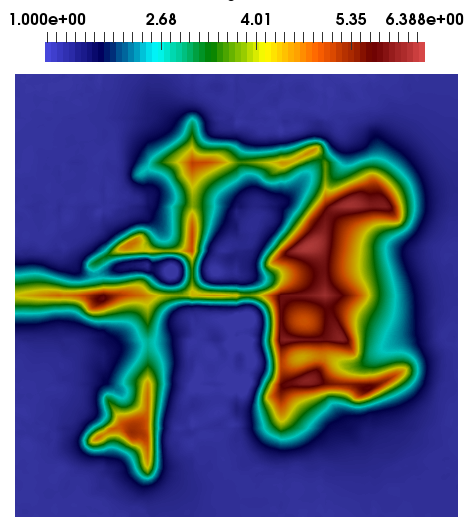}
\includegraphics[width=0.24\linewidth]{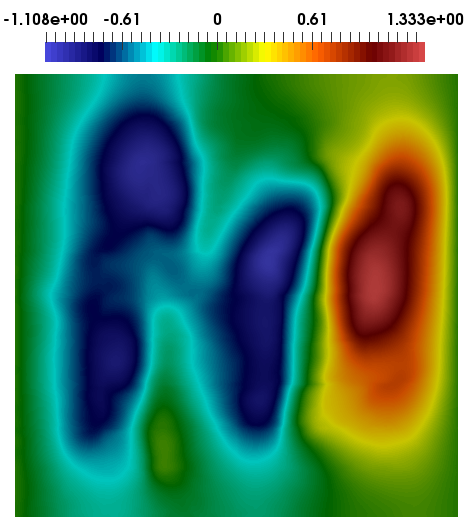}
\includegraphics[width=0.24\linewidth]{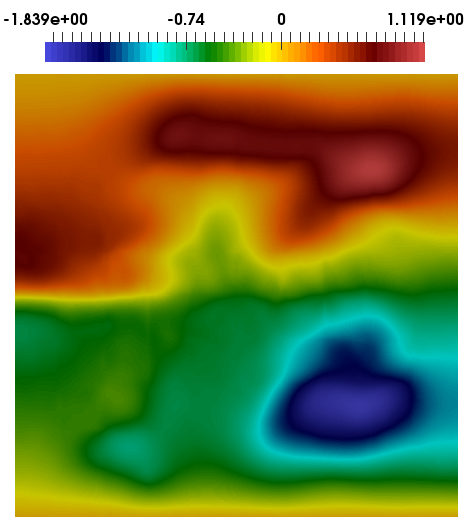} \\
\caption{Poroelasticity in multicontinuum media: Distribution of pressure for first and second continuum and displacement along $X$ and $Y$ at final time for fractional order derivative $\alpha = 0.9$ (from left to right).
First row: exact solution. Second row: multiscale solution(12 multiscale basis functions).}
\label{pic-5}
\end{figure}

\begin{figure}[h!]
\centering
\vspace{2.0 mm}
\includegraphics[width=0.24\linewidth]{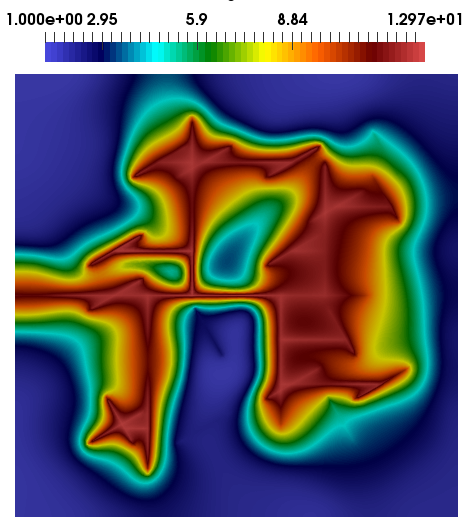}
\includegraphics[width=0.24\linewidth]{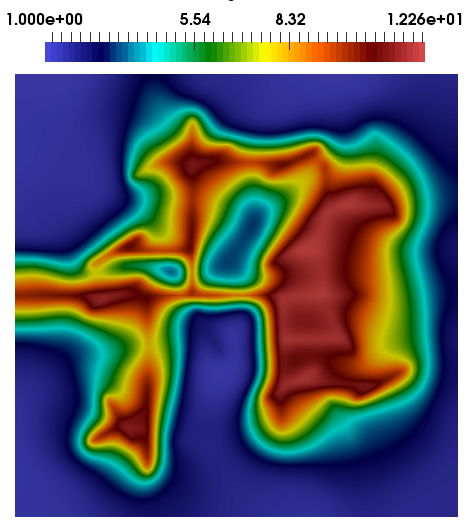}  
\includegraphics[width=0.24\linewidth]{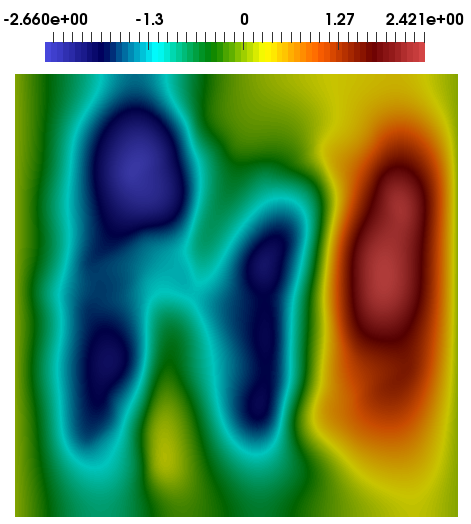}
\includegraphics[width=0.24\linewidth]{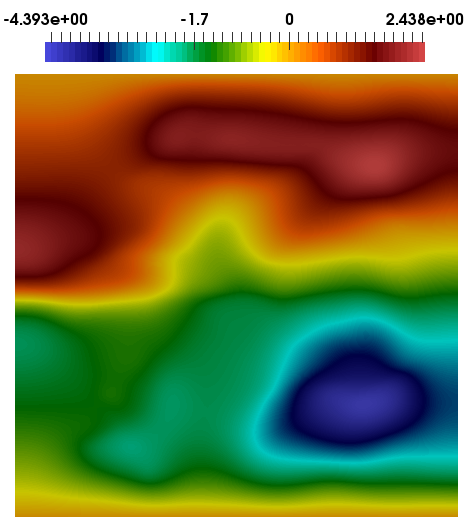} \\
\includegraphics[width=0.24\linewidth]{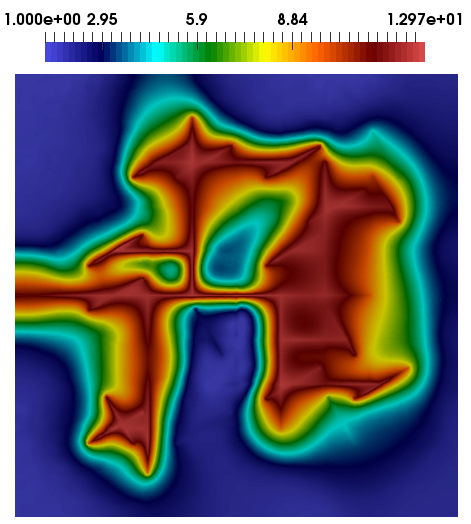}
\includegraphics[width=0.24\linewidth]{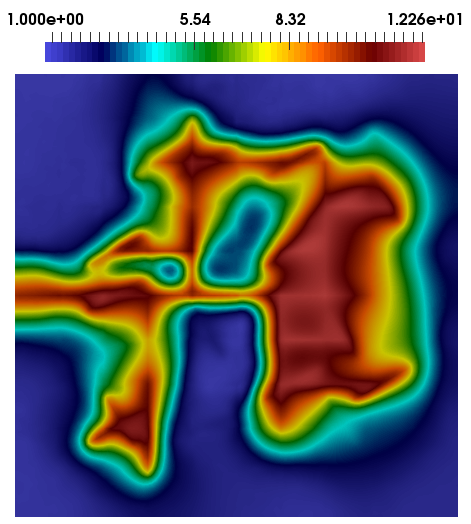}
\includegraphics[width=0.24\linewidth]{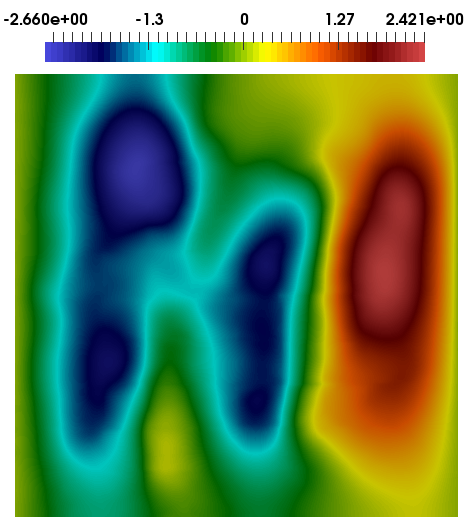}
\includegraphics[width=0.24\linewidth]{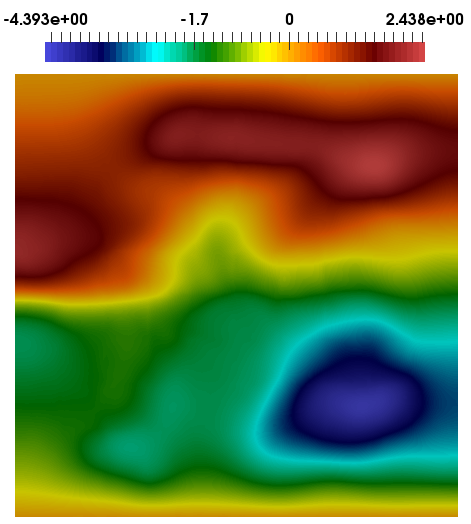} \\
\caption{Poroelasticity in multicontinuum media: Distribution of pressure for first and second continuum and displacement along $X$ and $Y$ at final time for fractional order derivative $\alpha = 1.0$ (from left to right).
First row: exact solution. Second row: multiscale solution(12 multiscale basis functions).}
\label{pic-6}
\end{figure}

The distribution of pressure for the first and second continua,
displacement along $x$ and $y$ directions at final time are presented
in Figures $\ref{pic-4}-\ref{pic-6}$. In the first row, we depict a
reference fine grid solution and multiscale solution with twelve
multiscale basis functions for the GMsFEM is presented in second
row. We observe good accuracy comparing the fine-scale solution with
the multiscale solution with twelve basis functions for displacement
along $x$ and $y$ direction and pressures for different fractional order derivatives. For the poroelasticity
problems in multicontinuum media, we also observe good convergences.

\begin{table}[h!]
\begin{center}
\begin{tabular}{ | c | c | c | c | c | c | c | c |}
\hline
$M$ & $DOF_H$ & $e_{L_2}^u$ (\%) & $e_{H_1}^u$  (\%) & $e_{L_2}^{p_1}$  (\%) & $e_{H_1}^{p_1}$  (\%) & $e_{L_2}^{p_2}$  (\%) & $e_{H_1}^{p_2}$  (\%) \\ \hline
1 & 484 & 26.201 & 64.703 & 17.893 & 83.608 & 13.191 & 78.062\\
2 & 968 & 14.965 & 43.514 & 13.705 & 68.893 & 9.745 & 65.278\\
4 & 1936 & 8.201 & 33.273 & 6.753 & 43.739 & 4.696 & 39.658\\
8 & 3872 & 4.142 & 23.807 & 3.350 & 27.219 & 2.315 & 26.418\\
12 & 5808 & 2.984 & 20.126 & 2.528 & 22.687 & 1.779 & 22.887 \\
16 & 7744 & 2.440 & 18.097 & 2.154 & 20.365 & 1.528 & 21.085\\
\hline
\end{tabular}
\end{center}
\caption{Relative errors for displacement and pressure with different numbers of multiscale basis functions. Poroelasticity in multicontinuum media with fractional order derivative $\alpha = 0.8$}
\label{tablemult1}
\end{table}

\begin{table}[h!]
\begin{center}
\begin{tabular}{ | c | c | c | c | c | c | c | c |}
\hline
$M$ & $DOF_H$ & $e_{L_2}^u$ (\%) & $e_{H_1}^u$  (\%) & $e_{L_2}^{p_1}$  (\%) & $e_{H_1}^{p_1}$  (\%) & $e_{L_2}^{p_2}$  (\%) & $e_{H_1}^{p_2}$  (\%) \\ \hline
1 & 484 & 22.431 & 53.276 & 15.457 & 73.825 & 14.153 & 66.143\\
2 & 968 & 11.122 & 31.456 & 10.879 & 58.526 & 9.623 & 53.550\\
4 & 1936 & 4.931 & 21.444 & 4.369 & 33.371 & 3.885 & 30.638\\
8 & 3872 & 2.264 & 14.665 & 1.866 & 19.661 & 1.678 & 18.949\\
12 & 5808 & 1.629 & 12.265 & 1.355 & 16.404 & 1.234 & 16.020 \\
16 & 7744 & 1.335 & 10.942 & 1.167 & 14.948 & 1.061 & 14.730\\
\hline
\end{tabular}
\end{center}
\caption{Relative errors for displacement and pressure with different numbers of multiscale basis functions. Poroelasticity in multicontinuum media with fractional order derivative $\alpha = 0.9$}
\label{tablemult2}
\end{table}

\begin{table}[h!]
\begin{center}
\begin{tabular}{ | c | c | c | c | c | c | c | c |}
\hline
$M$ & $DOF_H$ & $e_{L_2}^u$ (\%) & $e_{H_1}^u$  (\%) & $e_{L_2}^{p_1}$  (\%) & $e_{H_1}^{p_1}$  (\%) & $e_{L_2}^{p_2}$  (\%) & $e_{H_1}^{p_2}$  (\%) \\ \hline
1 & 484 & 20.226 & 44.876 & 12.573 & 63.923 & 12.272 & 56.038\\
2 & 968 & 8.981 & 23.888 & 7.720 & 48.952 & 7.391 & 43.954\\
4 & 1936 & 3.255 & 14.252 & 2.404 & 25.288 & 2.290 & 23.452\\
8 & 3872 & 1.440 & 9.519 & 0.922 & 14.563 & 0.904 & 13.362\\
12 & 5808 & 1.051 & 7.925 & 0.667 & 12.349 & 0.644 & 11.017 \\
16 & 7744 & 0.876 & 7.059 & 0.579 & 11.422 & 0.553 & 10.100\\
\hline
\end{tabular}
\end{center}
\caption{Relative errors for displacement and pressure with different numbers of multiscale basis functions. Poroelasticity in multicontinuum media with fractional order derivative $\alpha = 1.0$}
\label{tablemult3}
\end{table}

\section{Conclusion}

In this paper, a mathematical formulation is introduced  for
poroelasticity problems in fractured and heterogeneous media with the
time fractional derivatives. We assume the media have a multiscale
nature and develop a computational macroscale model.  A \rv{finite}
difference approximation of the Caputo fractional time derivative is
\rv{adopted} for flow and mechanics. Due to the time fractional order,
the resulting system has a memory and requires storing the solutions
at previous time steps. Because of multiple scales, we use the GMsFEM as a computational model. 
For the GMsFEM, one needs multiscale basis functions and a global formulation.
 We construct multiscale basis functions for the approximation of
 pressure and displacement and solve the problem on the coarse
 grid. The multiscale approach uses the Discrete Fracture Model to
 resolve the fractures on a fine grid. The numerical examples are
 presented to verify the efficiency of the proposed difference schemes
 for two-dimensional problem. We provide comparison results using
 different numbers of basis functions for the pressures in each
 continuum and the displacement between the multiscale and fine-scale solutions with different fractional order derivative. Our results show that the proposed method can give accurate solutions.

\section*{Acknowledgments}
The research of DS is supported in part by National Research Foundation (NRF-2017R1A2B3012506). The works of AA and AT are supported by North-Caucasus Center for Mathematical Research under agreement N. 075-02-2021-1749 with the Ministry of Science and Higher Education of the Russian Federation.
AT is supported by Russian government project Science and Univer-sities 121110900017-5 aimed at supporting junior laboratories. MV work is supported by the mega-grant of the Russian Federation Government №14.Y26.31.0013.

\bibliographystyle{abbrv}
\bibliography{lit}

\end{document}